%% file: main.tex
\documentclass{article}

\usepackage{color}
\usepackage{ulem}
\usepackage{siunitx}
\usepackage{outlines}
\usepackage{amsmath}
\usepackage{xcolor}
\usepackage{caption}
\usepackage{subcaption}
\usepackage{bm}
\usepackage{natbib}
\usepackage{mathtools}
\mathtoolsset{showonlyrefs}
\usepackage{wrapfig}
\usepackage{gensymb}
\usepackage{float}
\usepackage{algorithm}
\usepackage{algpseudocode}
\usepackage{todonotes}
\usepackage{placeins}
\usepackage{amssymb}
\usepackage{stmaryrd}
\usepackage[margin=3cm]{geometry}
\usepackage{authblk}
\usepackage[skip=10pt plus1pt, indent=40pt]{parskip}

\newcommand{\bA}{\mathbf{A}}
\newcommand{\bB}{\mathbf{B}}
\newcommand{\bC}{\mathbf{C}}
\newcommand{\bCo}{\bC_{\rm o}}
\newcommand{\bCb}{\bC_{\rm b}}
\newcommand{\bD}{\mathbf{D}}
\newcommand{\bF}{\mathbf{F}}
\newcommand{\bFp}{\bF_{p}}
\newcommand{\bG}{\mathbf{G}}
\newcommand{\bGt}{\bG^{\rm T}}

\newcommand{\bI}{\mathbf{I}}
\newcommand{\bK}{\mathbf{K}}
\newcommand{\bM}{\mathbf{M}}
\newcommand{\bMl}{\overline{\mathbf{M}}}
\newcommand{\bPa}{\mathbf{P}_{\rm a}}
\newcommand{\bR}{\mathbf{R}}
\newcommand{\bS}{\mathbf{S}}
\newcommand{\bU}{\mathbf{U}}
\newcommand{\bV}{\mathbf{V}}

\newcommand{\bGamma}{\bm{\Gamma}}
\newcommand{\bLambo}{\bm{\Lambda}_{\rm o}}
\newcommand{\bLambb}{\bm{\Lambda}_{\rm b}}
\newcommand{\bLambs}{\bm{\Lambda}_{\rm s}}
\newcommand{\bLamba}{\bm{\Lambda}_{\rm a}}

\newcommand{\bkappa}{\bm{\kappa}}

\newcommand{\bx}{\mathbf{x}}
\newcommand{\bdx}{\delta\mathbf{x}}
\newcommand{\bxa}{\mathbf{x}_{\rm a}}
\newcommand{\bdxa}{\delta \mathbf{x}_{\rm a}}
\newcommand{\bxb}{\mathbf{x_{\rm b}}}
\newcommand{\by}{\mathbf{y}_{\rm o}}

\newcommand{\bdya}{\delta \mathbf{y}_{\rm a}}
\newcommand{\bd}{\mathbf{d}}
\newcommand{\bxt}{\mathbf{x_{\rm t}}}
\newcommand{\bepsb}{\bm{\epsilon}_{\rm b}}
\newcommand{\bepso}{\bm{\epsilon}_{\rm o}}
\newcommand{\bepsa}{\bm{\epsilon}_{\rm a}}

\newcommand{\rhob}{\rho_{\rm b}}
\newcommand{\rhoo}{\rho_{\rm o}}
\newcommand{\mb}{m_{\rm b}}
\newcommand{\mo}{m_{\rm o}}
\newcommand{\sigmao}{\sigma_{\rm o}}
\newcommand{\sigmab}{\sigma_{\rm b}}
\newcommand{\sigmaa}{\sigma_{\rm a}}
\newcommand{\lambo}{\lambda_{\rm o}}
\newcommand{\lambb}{\lambda_{\rm b}}
\newcommand{\lambs}{\lambda_{\rm s}}
\newcommand{\lamba}{\lambda_{\rm a}}
\newcommand{\lambsk}{\lambs^{(k)}}
\newcommand{\lambok}{\lambo^{(k)}}
\newcommand{\lambbk}{\lambb^{(k)}}
\newcommand{\lambak}{\lamba^{(k)}}
\newcommand{\etaaa}{\eta_{\rm \widetilde{a}, a}^{(k)}}
\newcommand{\etaoo}{\eta_{\rm \widetilde{o}, o}^{(k)}}
\newcommand{\etaob}{\eta_{\rm o, b}^{(k)}}

\newcommand{\practical}[1]{\widetilde{#1}}

\newcommand{\PR}{\practical{\bR}}
\newcommand{\PK}{\practical{\bK}}
\newcommand{\PS}{\practical{\bS}}
\newcommand{\PLamba}{\practical{\bm{\Lambda}}_{\rm_a}}

\newcommand{\PLambs}{\practical{\bm{\Lambda}}_{\rm_s}}
\newcommand{\PPa}{\practical{\mathbf{P}}_{\rm a}}
\newcommand{\Pmo}{\practical{m}_{\rm o}}
\newcommand{\Prhoo}{\practical{\rho}_{\rm o}}
\newcommand{\Psigmao}{\practical{\sigma}_{\rm o}}
\newcommand{\Psigmaa}{\practical{\sigma}_{\rm a}}
\newcommand{\Plambsk}{\practical{\lambda}_{\rm s}^{(k)}}
\newcommand{\Plambak}{\practical{\lambda}_{\rm a}^{(k)}}

\newcommand{\Plambok}{\practical{\lambda}_{\rm o}^{(k)}}

\newcommand{\F}[1]{\widehat{#1}}
\newcommand{\fbdya}{\F{\delta \mathbf{y}}_{\rm a}}
\newcommand{\fepso}{\F{\bm{\epsilon}}_{\rm o}}
\newcommand{\fdyak}{\F{\delta y}_{\rm a}^{(k)}}
\newcommand{\fdk}{\F{d}^{(k)}}

   \newcommand{\E}[1]{\mathbb{E}\! \left[#1\right]} 

   \newcommand{\G}[1]{\mathcal{G}\! \left(#1\right)}

   \newcommand{\J}[1]{\mathcal{J}\! \left(#1\right)}

   \newcommand{\R}{\mathbb{R}}

\newcommand{\eq}[1]{\mbox{$#1$}}

\title{On the impact of observation error correlations in data assimilation, with application to along-track altimeter data}
\author{Olivier Goux\textsuperscript{th}}
\author[1]{Anthony T. Weaver}
\author[1]{Selime Gürol}
\author[2]{Oliver Guillet}
\author[1]{Andrea Piacentini}
\affil[1]{CECI UMR 5318, CERFACS and CNRS, Toulouse, France}
\affil[2]{CNRM UMR 3589, Météo-France and CNRS, Toulouse, France}
\date{}                     
\begin{document}

\maketitle

\begin{center}
This document is a pre-print version of an article currently (March 11, 2025) under review at the \textit{Quaterly Journal or the Royal Meteorological Society}
\end{center}

\begin{abstract}
Data assimilation involves estimating the state of a system by combining observations from various sources with a background estimate of the state. The weights given to the observations and background state depend on their specified error covariance matrices.  Observation errors are often assumed to be uncorrelated even though this assumption is inaccurate for many modern data-sets such as those from satellite observing systems. As methods allowing for a more realistic representation of observation-error correlations are emerging, our aim in this article is to provide insight on their expected impact in data assimilation. First, we use a simple idealised system to analyse the effect of observation-error correlations on the spectral characteristics of the solution. Next, we assess the relevance of these results in a more realistic setting in which simulated along-track (nadir) altimeter observations with correlated errors are assimilated in a global ocean model using a three-dimensional variational assimilation (3D-Var) method. Correlated observation errors are modelled in the 3D-Var system using a diffusion operator.

When the correlation length scale of observation error is small compared to that of background error, inflating the observation-error variances can mitigate most of the negative effects from neglecting the observation-error correlations. Accounting for observation-error correlations in this situation still outperforms variance inflation since it allows small-scale information in the observations to be more effectively extracted and does not affect the convergence of the minimization. Conversely, when the correlation length scale of observation error is large compared to that of background error, the effect of observation-error correlations cannot be properly approximated with variance inflation. However, the correlation model needs to be constructed carefully to ensure the minimization problem is adequately conditioned so that a robust solution can be obtained. Practical ways to achieve this are discussed.

\end{abstract}

\section{Introduction}
\label{sec_introduction}

The background- and observation-error covariance matrices (\eq{\bB} and \eq{\bR}, respectively) are fundamental components of data assimilation systems for the atmosphere and ocean. Much research has been devoted to the specification of \eq{\bB}, to the extent that quite sophisticated representations of \eq{\bB} have been developed for operational systems \citep[\textit{e.g.}, for a review see][]{Bannister_2008}. Comparatively, \eq{\bR} has received much less attention, especially with regard to the specification of spatial correlations (off-diagonal elements in \eq{\bf R}), which are generally neglected. Neglecting spatial observation-error correlations, however, is not a good assumption for many modern observation datasets. Multiple studies describe significant observation-error correlations in, for example, radiances from the Infrared Atmospheric Sounding Interferometer \mbox{\citep{Bormann_2010}} and the Spinning Enhanced Visible and InfraRed Imager \mbox{\citep{Waller_2016, Michel_2018}}, sea-level anomalies from satellite altimeters \citep{Storto_2019}, Doppler radial winds \citep{Waller_2016b}, and satellite-retrieved sea surface temperature \citep{Reid_2020}. Moreover, particularly large spatial error correlations are to be expected in a new-generation of high-resolution satellite observations, such as wide-swath altimeter data \citep{Fernandez_2018, King_2021} and total surface current velocity measurements \citep{Waters_2024}.
The benefits of accounting for correlated observation error, particularly for the analysis of small spatial scales, have been highlighted in several studies \citep{Healy_2005, Stewart_2013, Rainwater_2015, Pinnington_2016}.\par

Inter-channel error correlations in satellite radiance observations are now accounted for in multiple Numerical Weather Prediction systems \citep{Weston_2014, Bormann_2016,Campbell_2017, Geer_2019}, with a positive impact on forecast scores. Accounting for spatial correlations, however, has proved to be more challenging, since, unlike observations from different satellite instrument channels, observations at different spatial locations cannot be easily separated into manageable, independent subsets of observations for computing correlations. Though still in early development, practical methods are emerging to account for horizontal correlations in \eq{\bR}, particularly in the context of variational data assimilation.

In variational data assimilation, a cost function is minimised to find the state of the system that best fits the background state and observations in a weighted least-squares sense, where the weighting matrices are defined by \eq{\bB^{-1}} and \eq{\bR^{-1}}. The minimization strategy is typically based on an iterative Gauss-Newton algorithm \mbox{\citep{Gratton_2007}}, which in the data assimilation literature is usually referred to as an incremental algorithm \citep{Courtier_1994}. Central to the Gauss-Newton algorithm  is a so-called inner loop that involves minimising a quadratic cost function using a conjugate gradient (CG) method. An important feature of the algorithm is that matrices do not need to be defined explicitly; only operators representing  matrix-vector products are required. Furthermore, it is customary to employ \eq{\bB} as a preconditioning matrix for CG \mbox{\citep[\textit{e.g.},][]{Derber_1989, Gurol_2014}}. In a \eq{\bB}-preconditioned CG algorithm, an operator needs to be specified for \eq{\bB} (or a factored form of \eq{\bB}) as well as for the inverse covariance matrix, \eq{\bR^{-1}}. The fact that \eq{\bR^{-1}}, and not \eq{\bR}, is required by the minimization algorithm brings additional computational challenges when accounting for spatially correlated observation errors.\par

To compute a matrix-vector product with \eq{\bR^{-1}}, \citet{Hu_2024} propose the local Singular Value Decomposition (SVD) - Fast Multipole Method (FMM), which is a variant of the original SVD-FMM method of \mbox{\cite{Hu_2021}} but designed to be more computationally efficient and to require a reduced number of parallel processor communications. The main purpose of their method is to reduce the cost of computing the product of a vector with \eq{\bR^{-1}}, where they assume the latter to be available in explicit matrix form. To do so, the domain is partitioned into `boxes' containing subsets of observations. The block-diagonal components of \eq{\bR^{-1}} that correspond to observations in the same box are stored and applied directly to sub-vectors associated with those observations. The remaining block components of \eq{\bR^{-1}} are approximated, by using truncated SVD for blocks that are associated with observations within a prescribed `interaction' area, or by explicitly setting the blocks to zero if the observations are located beyond the interaction area. Crucially, the method relies on the availability of \eq{\bR^{-1}}, which in \mbox{\cite{Hu_2024}} is obtained by inverting a covariance matrix built explicitly from a parametric covariance function.\par

\mbox{\citet{Guillet_2019}} consider a different approach, which uses knowledge of \eq{\bR} to construct a meaningful \eq{\bR^{-1}} operator, but does not require \eq{\bR} itself to be built explicitly (and then to be inverted), which for large data-sets is challenging.  This is achieved by considering a specific set of covariance functions from the Mat\'ern class for which analytical expressions of the inverse covariance operator are known and straightforward to represent numerically as a product of sparse matrices. As discussed in \mbox{\citet{Guillet_2019}} and references therein, the inverse covariance operator associated with these functions can be interpreted as the inverse of an iterative diffusion process. The diffusion tensor and number of diffusion iterations are adjustable parameters that control the spatial-range and shape of the associated Mat\'ern covariance functions.
To use this approach to model \eq{\bR^{-1}}, the (inverse) diffusion operator is discretised on a triangular mesh with nodes at the observation locations. The Finite Element Method (FEM) with standard \eq{\mathcal{P}_1}~elements is used for spatially discretizing the inverse diffusion operator. Since each node on the mesh depends only on its direct neighbours, \eq{\bR^{-1}} is sparse and scalable with large volumes of data on memory-distributed computers. With knowledge of the underlying covariance function, the parameters of the diffusion model can be estimated by fitting the function to samples of observation-error covariances obtained, for example, from Desroziers diagnostics \mbox{\citep{Desroziers_2005}}.\par

An indirect way to account for spatially correlated observation errors is to augment the observation vector with spatial derivatives of the observations in addition to the observations themselves and to use a diagonal approximation for the error covariance matrix of the augmented observation vector
\citep{Brankart_2009, Chabot_2015, Ruggiero_2016, Yaremchuk_2018, Bedard_2019}. In a continuous framework, \mbox{\citet[Section~4.1]{Guillet_2019}} show that, for specific parameter settings, this is equivalent to the diffusion approach for modelling a non-diagonal \eq{\bR^{-1}} for the non-augmented observations.
 Assuming that an estimate of \eq{\bR} is available and can be stored in matrix form, \mbox{\citet{Ruggiero_2016}} and \mbox{\citet{Yaremchuk_2018}} describe an approach for estimating the weights (the inverse error variances) of the augmented observation-error covariance matrix that best fits the estimate of
\eq{\bR} in a squared Frobenius norm sense. Although the method does not require a direct implementation of a non-diagonal \eq{\bR}, it substantially increases the effective size of the observation vector (by almost a factor of five in the example of \cite{Ruggiero_2016}) and is primarily designed to handle observations with structured distributions for which simple finite-difference formulae can be used to compute spatial derivatives. Adapting the method to sparse and heterogeneously distributed observations requires more sophisticated numerical techniques such as those presented in \mbox{\citet[Section~4.2]{Guillet_2019}}.

The purpose of the current article is threefold. First, we aim to complement other studies by providing further insight on the expected impact of correlated observation error in data assimilation. We address this initially by studying the optimal analysis equations in an idealised, one-dimensional (1D) framework where \eq{\bB} and \eq{\bR} are described by circulant matrices and hence diagonalizable in a Fourier basis. By analysing the equations in Fourier space, we are able to highlight the scale-dependent characteristics of the optimal solution, focussing on the role of \eq{\bR}. We then use the same idealised framework to analyse the sensitivity of the solution to a misspecification of \eq{\bR}. Second, we assess the practical relevancy of these results experimentally using a global-ocean three-dimensional variational data assimilation (3D-Var) configuration of the NEMOVAR system \citep{Mogensen_2009}, where assimilated sea-level anomalies from (nadir) altimeters are taken to have correlated observation errors along the satellite track. While the configuration is realistic, we design the experiments in an idealised way to help with the interpretation of the results in relation with the theoretical analysis presented earlier. We also pay particular attention to the convergence properties of the (\eq{\bB}-preconditioned) CG minimization, which previous studies, involving mainly theoretical analysis, have shown to be highly sensitive to a non-diagonal \eq{\bR} \citep{Haben_2011, Tabeart_2018, Tabeart_2020, Goux_2024}. Third, with this work, we provide a practical demonstration of the method of \cite{Guillet_2019} in an operationally-based data assimilation system. We view this as an important first step towards applying the method to account for spatially correlated error in other data-sets such as wide-swath altimeter and satellite-retrieved sea-surface temperature.

The structure of the article is as follows. In Section~2, we present the data assimilation problem and the idealised theoretical framework for studying the impact on the analysis of both correctly and incorrectly specified  representations of \eq{\bR}. In Section~\ref{sec_experimental}, we describe the ocean data assimilation system, focussing on the method for accounting for along-track correlated errors in altimeter observations. We provide more details on the method and its technical implementation in Appendix~\ref{app_implementation}. In Section~3, we also present the idealised experimental framework that we use for assessing the impact of the correlated observation-error model. We discuss the results of the experiments in Section~4. Finally, we provide a summary and conclusions in Section~5.

\section{A simplified theoretical framework}
\label{sec_framework}

In this section, we first recall the basic equations describing the optimal variational analysis in a weakly non-linear system. We then present the idealised framework that allows us to interpret the analysis equations in Fourier space and to investigate the sensitivity of the analysis to the observation-error covariances.

\subsection{The variational analysis equations}
\label{sec_formulation}
\input{./problem_formulation.tex}

\subsection{A Fourier decomposition of the analysis equations}
\label{sec_theoretical}
\input{./theoretical_framework.tex}

\section{Data assimilation system and experimental framework}  \label{sec_experimental}
\input{./experimental_framework.tex}

\section{Numerical Results}
\label{sec_results}

We use both the theoretical framework of Section~\ref{sec_theoretical} and the experimental framework described in Section~\ref{sec_experimental} to interpret results for two general situations that involve different observation-error correlation length scales. First, in Section~\ref{sec_small_corr}, we consider a situation where the observation-error correlations have length scales that are short compared to those of the background. This is the expected situation with nadir altimeter data \citep[\textit{e.g.}, ][]{Storto_2019}. Second, in Section~\ref{sec_large_corr}, we consider a situation where  the observation-error correlations have length scales that are large compared to those of the background. This situation is to be expected with wide-swath altimeter data such as SWOT (Surface Water Ocean Topography) due to observation errors being significantly correlated at distances up to 1000~km \citep{Fernandez_2018}. In Section~\ref{sec_inflation}, we show how the presence of different physical variables, which was not accounted for in the theoretical framework, is one of the factors limiting the effectiveness of variance inflation compared to the use of a non-diagonal \eq{\PR}. In both situations, we choose the correlation model smoothness parameters such that \eq{\mo < \mb}. This choice is motivated by results from \citet{Goux_2024} who showed an extreme sensitivity of the condition number and, as a consequence, the minimization convergence when \eq{\mo \geq\mb}. The sensitivity to \eq{\mo} is discussed in Section~\ref{sec_recond}.

\subsection{Case~1. Short observation-error correlation length scales: \eq{\rhoo < \rhob} }
\label{sec_small_corr}
\input{./small_corr.tex}

\subsection{Case~2. Large observation-error correlation length scales: \eq{\rhoo > \rhob}}
\label{sec_large_corr}
\input{./large_corr.tex}

\subsection{Optimal variance inflation}
\label{sec_inflation}
\input{./app_inflation.tex}

\subsection{Sensitivity of \eq{\PR^{-1}}}
\label{sec_recond}
\input{./stability.tex}

\FloatBarrier
\section{Summary and conclusions}
\label{sec_conclusion}
\input{./conclusion.tex}

\section*{acknowledgements}
This work has benefitted from funding from the Copernicus Climate Change Service (C3S), one of six services of the European Union’s Copernicus Programme, which is implemented by ECMWF on behalf of the European Commission. It is a contribution to the contract C3S\_601\_INRIA: ``Advancing ocean data assimilation methodology for climate applications''. Additional support has been provided by the French National Programme LEFE/INSU (Les Enveloppes Fluides et l’Environnement/Institut National des Sciences de l'Univers) through the scientific themes MANU (M\'ethodes Math\'ematiques et NUm\'eriques) and GMMC (Groupe Mission Mercator/Coriolis).

\section*{conflict of interest}
This study does not have any conflicts to disclose.

\section*{Data availability}
Data sharing is not applicable to this article as no new data were created or analysed in this study.

\bibliographystyle{apalike}
\bibliography{./references}

\newpage

\begin{appendix}

\section{Implementation of the finite-element method}
\FloatBarrier
\label{app_implementation}
\input{./app_implementation.tex}

\end{appendix}

\end{document}

%% file: problem_formulation.tex
Data assimilation methods seek to combine a background state \eq{\bxb\in \mathbb{R}^n} and observations \eq{\by \in \mathbb{R}^p} in a statistically optimal manner \mbox{\citep{Fletcher_2022}}. We assume the background state to be the sum of the true state \eq{\bxt\in \mathbb{R}^n} that we are attempting to estimate, and an additive Gaussian error \eq{\bepsb\in \mathbb{R}^n} with mean zero and covariance matrix \eq{\bB\in \mathbb{R}^{n \times n}}:
\begin{equation}
   \bxb = \bxt + \bepsb.
   \label{eq_def_xb}
\end{equation}
We define a generalised observation operator \eq{\mathcal{G}} such that, for any estimate of the state of the system \eq{\bx \in \mathbb{R}^{n}}, the vector \eq{\G{\bx}} represents how the state of the system would appear through the observation network. In 4D-Var, \eq{\mathcal{G}} will include the time-dependent model operator as well as the observation operator itself \citep{Courtier_1994}.  We assume the observations to be the sum of the counterpart of the true state in observation space and an additive Gaussian error \eq{\bepso\in \mathbb{R}^p} with mean zero and covariance matrix \eq{\bR \in \mathbb{R}^{p \times p}}:
\begin{equation}
   \by = \G{\bxt} + \bepso.
   \label{eq_def_yo}
\end{equation}

In variational data assimilation, the analysis is obtained by minimizing a cost function \eq{\mathcal{J}} that measures the simultaneous fit of \eq{\bx} to \eq{\bxb} and \eq{\by}: 
\begin{equation}
   \J{\bx} = \frac{1}{2} \|\bx - \bxb\|_{\bB^{-1}}^2  + \frac{1}{2} \|\G{\bx} - \by \|_{\bR^{-1}}^2,
\end{equation}
which can be written in terms of an increment \eq{\bdx = \bx-\bxb} as
\begin{equation}
\J{\bxb +\bdx} = \frac{1}{2} \|\bdx\|_{\bB^{-1}}^2  + \frac{1}{2} \|\G{\bxb+\bdx} - \by \|_{\bR^{-1}}^2.
   \label{eq:nonquadJ}
\end{equation}
The inverse error covariance matrices \eq{\bB^{-1}} and \eq{\bR^{-1}} define symmetric, positive-definite (SPD) weighting matrices for the squared norm $\| \bx \|_{\bA^{-1}}^2 = \bx^{\rm T} \bA^{-1} \bx$.
The minimizing solution of this cost function with nonlinear constraints is generally obtained by minimizing successive quadratic cost functions involving linearised constraints \citep{Courtier_1994}. 
 For the sake of simplicity, we will focus on the first of these approximating quadratic cost functions, which is obtained by replacing \eq{\mathcal{G}} with a linear operator
\begin{equation}
\G{\bxb + \bdx} \simeq \G{\bxb} + \bG\bdx
\label{eq:Gapprox}
\end{equation}
where \eq{\bG} is the Jacobian matrix of \eq{\mathcal{G}} (or of a simplified operator that approximates  \eq{\mathcal{G}}) evaluated at \eq{\bx = \bxb}.
Substituting Equation~\eqref{eq:Gapprox} in Equation~\eqref{eq:nonquadJ} yields the quadratic cost function
\begin{align}
   J(\bdx) &= \frac{1}{2} \|\bdx\|_{\bB^{-1}}^2  + \frac{1}{2} \|\G{\bxb}+\bG\bdx - \by \|_{\bR^{-1}}^2 \\
   &= \frac{1}{2} \|\bdx\|_{\bB^{-1}}^2  + \frac{1}{2} \|\bG\bdx - \bd \|_{\bR^{-1}}^2,
\end{align}
where
\begin{equation}
\eq{\bd = \by - \G{\bxb}}
\label{eq_def_innovation}
\end{equation}
is the innovation vector.
Under these assumptions, the expression for the analysis increment, \eq{\bdxa=\bxa -\bxb}, that minimizes the cost function \eq{J} (associated with the first outer loop) is well known \citep[\textit{e.g.},][]{Fletcher_2022}:
\begin{equation}
   \bdxa = \bK\bd
\end{equation}
where
\begin{equation}
   \bK = \bB\bGt\left(\bG\bB\bGt + \bR \right)^{-1} = \bPa\bGt\bR^{-1},
 \label{eq:K}
\end{equation}
is the Kalman gain matrix, and
\begin{equation}
\bPa =\left(\bB^{-1} + \bGt\bR^{-1}\bG \right)^{-1}
\label{eq_Pa_opt}
\end{equation}
is the analysis-error covariance matrix (in the linear approximation). In observation space, the analysis increment is given by
\begin{equation}
   \bdya = \bG \bdxa = \bS \bd \label{eq_BLUE}
\end{equation}
where 
\begin{equation}
   \bS  = \bG\bK = \bG\bB\bGt\left(\bG\bB\bGt + \bR \right)^{-1} = \bG \bPa \bGt \bR^{-1}
   \label{eq_def_S}.
\end{equation}
Equation~\eqref{eq_def_S} is called the \textit{sensitivity} matrix since \eq{\bS = \partial \bdya/\partial \bd} and hence its elements quantify the sensitivity of the analysis increment (in observation space) with respect to the innovations \citep{Cardinali_2004, Fowler_2018}.\par

%% file: theoretical_framework.tex
To gain insight on the effect of the Kalman gain with a non-diagonal $\mathbf{R}$, we consider a simplified framework where instructive theoretical results can be derived. In particular, we assume the following:
\begin{enumerate}
   \item the domain is 1D and periodic;
   \item the state space consists of a single variable defined at discrete locations on a uniform grid;
   \item observations of the state variable are direct and available at uniformly distributed locations (\textit{i.e.}, \eq{\bG} is a selection operator);
   \item the background- and observation-error covariances are homogeneous.
\end{enumerate}
From the last assumption, we can define
\begin{align}
   \bR &= \sigmao^2 \bCo \\
\mbox{and} \hspace{0.5cm}  \bG\bB\bGt &= \sigmab^2 \bCb,
\end{align}
where \eq{\sigmao^2} is the observation-error variance, \eq{\sigmab^2} is background-error variance in observation space, \eq{\bCo} is the observation-error correlation matrix, and \eq{\bCb} is the background-error correlation matrix in observation space. The main interest of this idealised framework is that \eq{\bB}, \eq{\bR} and \eq{\bG\bB\bGt}, and their associated correlation matrices, are circulant matrices, meaning that each of their rows can be obtained as a cyclic permutation of the previous row \citep[see Lemma 1 in][]{Goux_2024}. \par

All circulant matrices are diagonal in a Fourier basis \citep{Gray_2005}. Let \eq{\bFp} be the \mbox{\eq{p} $\times$ \eq{p}} matrix whose element on the \eq{k}-th row and \eq{\ell}-th column is defined by the Fourier basis function
\begin{equation}
   [\bFp]_{k,\ell} = \frac{1}{\sqrt{p}}e^{\, j 2\pi \frac{k\ell}{p}}, \label{eq_def_Fm}
\end{equation}
where \eq{j} denotes the unit imaginary number (\eq{j^2=-1}). We can define the Fourier transform of any vector \eq{\bx} as
\begin{equation}
  \F{\bx}= \bFp^{\rm H} \;\bx,
\end{equation}
where the superscript `H' denotes conjugate (Hermitian) transpose. The matrix \eq{\bFp} is unitary; \textit{i.e.}, \eq{\bFp^{\rm H} = \bFp^{-1}}, so that \eq{\bFp^{\rm H}\bFp = \bFp\bFp^{\rm H}  = \bI}. \par

The matrices \eq{\bR} and \eq{\bG\bB\bGt} are both diagonal in the Fourier basis defined by \eq{\bFp}, and can thus be expressed as
\begin{align}
   \bR &= \bFp\bLambo\bFp^{\rm H}\\
\mbox{and} \hspace{0.5cm}    \bG\bB\bG^{\rm T} &= \bFp\bLambb\bFp^{\rm H},
\end{align}
where the matrices \eq{\bLambo} and \eq{\bLambb} are diagonal. The \eq{k}-th diagonal element of \eq{\bLambo} is \eq{\sigmao^2} \eq{\lambok} where \eq{\lambok} is the \eq{k}-th eigenvalue of \eq{\bCo}. Similarly, the \eq{k}-th diagonal element of \eq{\bLambb} is \eq{\sigmab^2} \eq{\lambbk} where \eq{\lambbk} is the \eq{k}-th eigenvalue of \eq{\bCb}.\par

Each eigenvector of \eq{\bR} (\textit{i.e.}, each column of \eq{\bFp}) can be associated with a different wavenumber \eq{p/k} as seen in Equation \eqref{eq_def_Fm}. To facilitate physical interpretation, this non-dimensional wavenumber can be converted into a `spatial scale' by multiplying it by \eq{h}, the distance between consecutive observation locations. This spatial scale \eq{ph/k} will be used as the horizontal axis in figures showing eigenvalue spectrums of circulant matrices throughout this article. An eigenvalue of \eq{\bR} (\textit{i.e.}, diagonal element of \eq{\bLambo}) can be interpreted as the {\it scale-dependent} observation-error variance at the spatial scale defined by the associated eigenvector. Specifically, \eq{\bLambo} corresponds to the covariance matrix of \eq{\fepso}, the Fourier transform of the observation error:
\begin{equation}
    \bLambo \, = \, \bFp^{\rm H}\bR \bFp
    \, = \, \bFp^{\rm H}\E{\bepso\bepso^{\rm T}}\bFp
    \, = \, \E{\bFp^{\rm H}\bepso\bepso^{\rm T} \bFp}
    \, = \, \E{\fepso \fepso^{\, \rm H}},
\end{equation}
where \eq{\E{}} denotes the mathematical expectation operator. Similarly, the eigenvalues of \eq{\bG\bB\bG^{\rm T}} describe the scale-dependent variances of \eq{\bG\bepsb} (\textit{i.e.}, the variances of the Fourier transform of the background error mapped to observation space):
\begin{equation}
    \bLambb \, = \, \bFp^{\rm H}\bG\bB\bGt \bFp
    \, = \, \bFp^{\rm H}\bG\E{\bepsb\bepsb^{\rm T}}\bGt\bFp
    \, = \, \E{(\bFp^{\rm H}\bG\bepsb)(\bFp^{\rm H}\bG\bepsb)^{\rm H}}.
\end{equation}

As \eq{\bG\bB\bG^{\rm T}} and \eq{\bR} are both circulant, the sensitivity matrix \eq{\bS} is also circulant\footnote{The sum, product and inverse of circulant matrices are also circulant matrices.}, and thus diagonal in the Fourier basis defined by \eq{\bFp}:
\begin{equation}
   \bS= \bFp \bLambs \bFp^{\rm H},
   \label{eq_diag_S}
\end{equation}
where \eq{\bLambs} is a diagonal matrix. We can obtain an explicit expression for the eigenvalues of \eq{\bS} by combining Equations~\eqref{eq_def_S} and \eqref{eq_diag_S} and using the orthogonality property of \eq{\bFp} to yield
\begin{equation}
   \bLambs = \bLambb \left( \bLambb + \bLambo \right)^{-1} = \left( \bLambb^{-1} + \bLambo^{-1} \right)^{-1} \bLambo^{-1}.
 \label{eq_Lambdas}
\end{equation}
 The \eq{k}-th diagonal element of \eq{\bLambs}, which corresponds to the \eq{k}-th eigenvalue of \eq{\bS}, can then be expressed as
\begin{equation}
   \lambsk =  \frac{1}{1+\nu^{(k)}} \quad\text{where}\quad \nu^{(k)} = \frac{\sigmao^2\lambok}{\sigmab^2\lambbk}.
   \label{eq:lambda_sk}
\end{equation}

 Using Equation~\eqref{eq_diag_S}, the Fourier transform of Equation~\eqref{eq_BLUE} has the simple form
\begin{align}
   \fbdya &= \bFp^{\rm H}\bdya = \bFp^{\rm H} \;\bS \;\bd = \bFp^{\rm H} \;\bS \;\bFp\bFp^{\rm H} \;\bd = \bLambs \F{\bd}. \label{eq_fdya}
\end{align}
 Let \eq{\fdyak} and \eq{\widehat{d}^{(k)}} denote the \eq{k}-th elements of the vectors \eq{\fbdya} and \eq{\F{\bd}}, respectively.
The relationship between the analysis increment in observation space and the innovation is trivial in Fourier space:
\begin{equation}
  \fdyak = \lambsk \widehat{d}^{(k)} = \left( \frac{1}{1+\nu^{(k)}} \right) \fdk. \label{eq_spectral_blue}
\end{equation}

Note that Equation~\eqref{eq_spectral_blue} has the same basic form as the expression of the analysis increment of a scalar quantity, obtained by combining a background estimate and a single direct observation \citep[\textit{e.g.},][]{Fisher_2007}:
\begin{equation}
   \delta y_{\rm a} = \delta x_{\rm a} = \left( \frac{1}{1+ \nu} \right) d
   \label{eq_scalar_blue}, \quad \text{with}\quad \nu = \frac{\sigmao^2}{\sigmab^2}.
\end{equation}
In the scalar case, \eq{\sigmao^2} and \eq{\sigmab^2} fully describe the observation- and background-error variance, respectively. In our case, however, the error correlations determine effective error variances  at each spatial scale through the eigenvalues \eq{\lambok} and \eq{\lambbk}. In turn, this determines a different weighting between the background and observations at each spatial scale. The eigenvalue \eq{\lambsk} can be interpreted as the \textit{sensitivity} of the analysis to the observations at the \eq{k}-th spatial scale. Hereafter, we will refer to \eq{\big( \lambsk \big)_{k \in \llbracket 1, p \rrbracket}} as the scale-dependent sensitivities. \par

The covariance matrix of the analysis error in observation space is also circulant and diagonalised by \eq{\bFp}:
\begin{equation}
   \bG\bPa\bGt = \bFp\bLamba\bFp^{\rm H}, \label{eq_def_Lamba}
\end{equation}
where \eq{\bLamba} is a diagonal matrix. Its \eq{k}-th element is denoted \eq{\sigmaa^2\lambak} and can be interpreted as a variance at a specific spatial scale of the analysis error in observation space. From Equations~\eqref{eq_def_S}, \eqref{eq_diag_S} and \eqref{eq_Lambdas}, it follows that
\begin{align}
    \bLamba = \left( \bLambb^{-1} + \bLambo^{-1} \right)^{-1}.\label{eq_Lamba_opt}
\end{align}
The \eq{k}-th diagonal element of \eq{\bLamba} is thus given by
\begin{align}
   \sigmaa^2\lambak =  \left(\frac{1}{\sigmab^2\lambbk} + \frac{1}{\sigmao^2\lambok}\right)^{-1}.
   \label{eq_analysis_error_spectral}
\end{align}

\begin{figure}[h]
\begin{center}
   \includegraphics[width=0.5\textwidth]{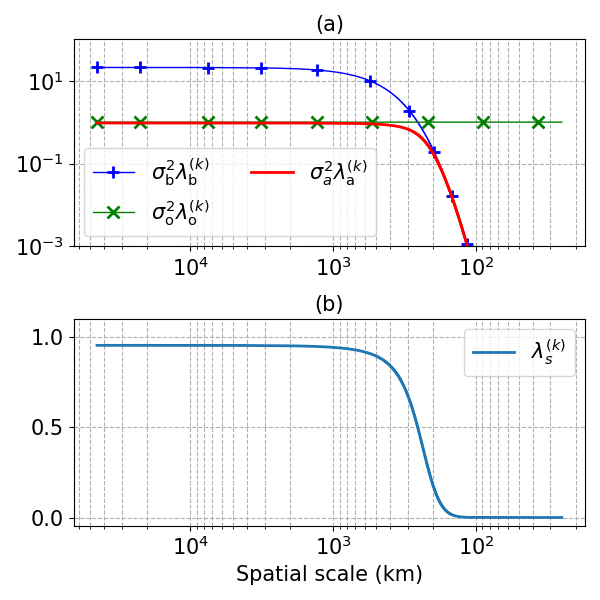}
      \caption{(a) Scale-dependent variance of the background error (\eq{\sigmab^2\lambbk}), observation error (\eq{\sigmao^2\lambok}), and  analysis error (\eq{\sigmaa^2\lambak}; Equation~\eqref{eq_analysis_error_spectral}), as a function of spatial scale.  (b) Scale-dependent sensitivity (\eq{\lambsk}; Equation~\eqref{eq:lambda_sk}) as a function of spatial scale. In this idealised 1D problem, the domain size is 40000~km; the background-error correlation function is approximately Gaussian with a length scale of 200~km; the observation errors are uncorrelated; and the background- and observation-error variances are both equal to one (see text for further description). The spatial scale is computed as \eq{ph/k}; \textit{i.e.}, the inverse of the non-dimensional wavenumber \eq{k/p} from Equation~\eqref{eq_def_Fm} multiplied by the distance between observations (\eq{h = 25}~km).
      The scale on the vertical axis has been cut-off at \eq{10^{-3}} (the minimum background-error variance reaches \eq{10^{-13}} for the smallest resolved spatial scale). }
   \label{fig_scale_decomp_1}
\end{center}
\end{figure}
Equation~\eqref{eq_analysis_error_spectral} is similar to the expression for the analysis-error variance for the scalar problem described earlier:
\begin{equation}
   \sigmaa^2 \, = \, \left(\frac{1}{\sigmab^2} + \frac{1}{\sigmao^2}\right)^{-1}.
   \label{eq_analysis_error_scalar}
\end{equation}
The equation pairs [\eqref{eq_scalar_blue}, \eqref{eq_analysis_error_scalar}] and [\eqref{eq_spectral_blue}, \eqref{eq_analysis_error_spectral}] highlight the similarity between the scalar example and what occurs at each spatial scale in the circulant framework.\par

We illustrate the properties of Equations~\eqref{eq_spectral_blue} and \eqref{eq_analysis_error_spectral} with a specific example. Consider a 1D domain of length 40000~km for which there is an observation every 25~km. The domain can be interpreted as a latitude circle at the Equator where observations are available every quarter of a degree in longitude. Both the observation- and background-error variances are set to one unit (\eq{\sigmao^2 = \sigmab^2 = 1}). The background-error correlation function is taken to be a 10-th order auto-regressive (quasi-Gaussian) function with a length scale of 250~km (see Section~\ref{sec_implementation}). The eigenvalues \eq{\lambbk} of the associated correlation matrix are computed using the analytical expression given by Equation~(16) in \citet{Goux_2024}. The observation error is taken to be uncorrelated, implying that \eq{\lambok = 1} for all $k$.

In this idealised example, the Kalman gain -- or more precisely the sensitivity matrix -- acts as a spectral filter to determine how damped the innovation should be at each spatial scale in order to produce an analysis increment with minimum error variance. As such, the analysis contains the best compromise between background and observations {\it at each spatial scale}. Figure~\ref{fig_scale_decomp_1}(a) shows the scale-dependent variances of observation error (\eq{\sigmao^2\lambok}), background error (\eq{\sigmab^2\lambbk}), and analysis error (\eq{\sigma^2\lambak}) as a function of spatial scale, while Figure~\ref{fig_scale_decomp_1}(b) shows the scale-dependent sensitivity of the analysis to the observations (\eq{\lambsk}) as a function of spatial scale. From these figures, we can make the following basic observations.
At {\it large} (\eq{>10^3}~km) spatial scales,
\eq{\sigmaa^2\lambak \simeq \sigmao^2 \lambok \ll \sigmab^2\lambbk} and \eq{\lambsk \simeq 1}, so that \eq{\fdyak \simeq \fdk} from Equation~\eqref{eq_spectral_blue}, and thus the analysis at observation locations fits the observations.
At {\it small} (\eq{<150}~km) spatial scales, \eq{\sigmaa^2\lambak \simeq \sigmab^2 \lambbk \ll \sigmao^2\lambok} and \eq{\lambsk \simeq 0}, so that \eq{\fdyak \simeq 0}, and thus the analysis at observation locations remains close to the background. At  spatial scale 250~km (\eq{k=180}), where the background- and observation-error variance intersect,
\eq{\sigmaa^2\lambak = \frac{1}{2} \sigmab^2 \lambbk = \frac{1}{2}\sigmao^2\lambok} and
\eq{\lambsk=1/2}, so that \eq{\fdyak = \frac{1}{2}\fdk}, and thus the analysis at observation locations is half-way between the background and observations.
\par


Equations~\eqref{eq_spectral_blue} and \eqref{eq_analysis_error_spectral} assume  perfect knowledge of the background- and observation-error covariances as well as perfect practical capability to use these covariances to define the Kalman gain. In practice, the error covariance matrices applied in the Kalman gain (which is defined implicitly in variational data assimilation through the iterative minimization of the cost function) are not exactly representative of the (unknown) `true' observation- and background-error covariance matrices, \eq{\bR} and \eq{\bB}. In the next section, we revisit the analysis above but considering the case where the observation-error covariances are misspecified. 

\subsection{Misspecified observation-error covariances}
\label{sec_misspec}
In this section, we study the effect of misspecifying the observation-error covariances, drawing inspiration from \citet{Fisher_2007}, who studied the effect of misspecifying the background-error variance in a scalar example. We introduce the notation \eq{\PR} to denote the {\it practical} observation-error covariance matrix used in the gain matrix, which can potentially differ from the {\it true} observation-error covariance matrix, \eq{\bR}. When using \eq{\PR} instead of \eq{\bR}, the resulting gain, sensitivity and analysis-error covariance matrices are suboptimal, and will be denoted by \eq{\PK}, \eq{\PS}, and \eq{\PPa}, respectively, to distinguish them from their `true' counterparts \eq{\bK}, \eq{\bS}, and \eq{\bPa} that are defined in terms of \eq{\bR}. Eigenvalues and variances associated with the suboptimal (practical) matrices \eq{\PR}, \eq{\PS} and \eq{\PPa} are also denoted with a tilde. As the focus here is on observation error, we assume that the practical and true background-error covariance matrices are equal.\par

With potentially misspecified observation-error covariances, we can define the practical gain matrix, the sensitivity matrix and the eigenvalues of the latter using expressions analogous to those in the optimal case:
\begin{align}
   \PK = \bB\bGt\left( \bG\bB\bGt + \PR\right)^{-1}, \quad \PS = \bG \PK
\end{align}
and
\begin{align}
   \Plambsk = \left(1+ \frac{\Psigmao^2\Plambok}{\sigmab^2\lambbk}\right)^{-1}.
   \label{eq_sensitivity_spectral_subopt}
\end{align}
The expression for the practical analysis-error covariance matrix, however, becomes more complex as it involves both the practical and true error covariance matrices. For any \eq{\PK}, potentially different from the optimal \eq{\bK}, the practical analysis-error covariance matrix \eq{\PPa} can be defined as \citep[\textit{e.g.}, Equation~4.2.12 in ][]{Gelb_1974}
\begin{equation}
    \PPa = (\bI -\PK\bG)\bB(\bI - \PK\bG)^{\rm T} + \PK \, \bR \, \PK^{\rm T}. \label{eq_PPa}
\end{equation}
As in the optimal case, we focus on the analysis-error covariance matrix in observation space, which is given by
\begin{equation}
    \bG\PPa\bGt = (\bI -\PS)\bG\bB\bGt(\bI - \PS) + \PS \, \bR \, \PS^{\, \rm T}.
\end{equation}
As in the optimal case (Equation~\eqref{eq_def_Lamba}), we can use \eq{\bFp} to diagonalize the left-hand side of Equation~\eqref{eq_PPa}:
\begin{equation}
  \bFp \PLamba \bFp^{\rm H} = \bG\PPa\bGt \label{eq_PPa_PLamba}
\end{equation}
where \eq{\PLamba} is a diagonal matrix.
Combining Equations \eqref{eq_PPa} and \eqref{eq_PPa_PLamba}, and using the fact that \eq{\bFp^{\rm H}\bFp=\bI}, yields
\begin{equation}
    \PLamba = (\bI -\PLambs)\bLambb(\bI - \PLambs) + \PLambs\bLambo\PLambs. \label{eq_Lamba_subopt}
\end{equation}
From the expression for the elements of \eq{\PLambs} in Equation~\eqref{eq_sensitivity_spectral_subopt}, the diagonal elements of \eq{\PLamba} can be expressed as
\begin{equation}
   \Psigmaa^2\Plambak = \frac{\sigmab^2\lambbk \left( \Psigmao^2\Plambok\right)^{2} + \sigmao^2\lambok \left( \sigmab^2\lambbk\right)^{2}}{\left(\Psigmao^2\Plambok + \sigmab^2\lambbk\right)^{2}}.
   \label{eq_analysis_error_spectral_subopt}
\end{equation}

The elements \eq{\Psigmaa^2\Plambak} are the scale-dependent analysis-error variances resulting from the use of a suboptimal representation of the observation-error covariances.
To help interpret Equation~\eqref{eq_analysis_error_spectral_subopt}, we first normalize it by \eq{\sigmaa^2\lambak} in Equation~\eqref{eq_analysis_error_spectral}. This yields the expression
\begin{equation}
    \etaaa = \frac{1+ \left( \etaob \right)^2 \left( \etaoo \right)^2 + \left[ 1 + \left( \etaoo \right)^2 \right] \etaob} {1+ \left( \etaob \right)^2 \left( \etaoo \right)^2 +2\etaoo\etaob}
    \label{eq:etaaa}
\end{equation}
where \eq{\etaaa = \Psigmaa^2\Plambak / \sigmaa^2\lambak}, \eq{\etaoo = \Psigmao^2\Plambok / \sigmao^2\lambok} and \eq{\etaob = \sigmao^2\lambok / \sigmab^2\lambbk}.
It is easy to see from Equation~\eqref{eq:etaaa} that, as expected, \eq{\Psigmaa^2\Plambak \geq \sigmaa^2\lambak} (\eq{\etaaa \ge 1})  and \eq{\Psigmaa^2\Plambak = \sigmaa^2\lambak} if \eq{\etaoo=1}.\par

\begin{figure}[h]
\begin{center}
   \includegraphics[width=0.5\textwidth]{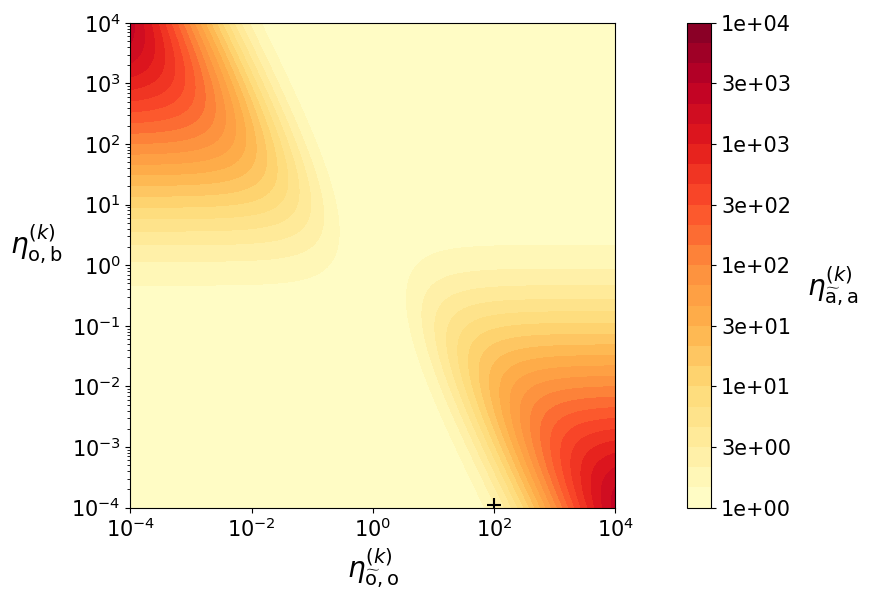}
   \caption{Ratio of the practical scale-dependent analysis-error variance to the optimal scale-dependent analysis-error variance \mbox{(\eq{\etaaa}; see Equation~\eqref{eq:etaaa})} where the former is associated with a misspecification of the observation-error variance at a given spatial scale (\eq{\etaoo\neq 1}). The black cross marks the point where \eq{\etaob = 10^{-4}}  and \eq{\etaoo=10^2} (referenced later in the main text).}
   \label{fig_error_sensitivity}
   \end{center}
\end{figure}

Figure~\ref{fig_error_sensitivity} shows \eq{\etaaa} plotted as a function of \eq{\etaoo} (horizontal axis) and \eq{\etaob} (vertical axis). Values of \eq{\etaoo} less (greater) than one correspond to an underestimation (overestimation) of the  observation-error variance at spatial scale $k$, while values of \eq{\etaob} less (greater) than one correspond to a more (less) accurate background than observation at spatial scale $k$.
The most striking feature of Figure~\ref{fig_error_sensitivity} is the asymmetry in magnitude between the top left/bottom right corners and the top right/bottom left corners. If the observations are more accurate than the background at a given spatial scale, overestimating \eq{\Psigmao^2\Plambok} is worse (in terms of increasing the analysis-error variance) than underestimating it; \textit{i.e.,} underfitting the observations is worse than overfitting them. If the observations are less accurate than the background at a given spatial scale, overestimating \eq{\Psigmao^2\Plambok} is better than underestimating it; \textit{i.e.,} underfitting the observations is better than overfitting them.\par

In other words, at any spatial scale, it is better to be too close to the most accurate source of information than too far from it.
Note that this last conclusion would also stand if the scale-dependent background-error variance was misspecified instead of the scale-dependent observation-error variance. This can be inferred from Figure~6 in \mbox{\citet{Fisher_2007}} which, although derived for a misspecification of the background-error variance in the scalar example, also translates to the scale-dependent background-error variance in our idealised circulant framework.\par

A less visually impressive but perhaps more impactful feature is that the iso-contours of the normalised analysis-error variance recede close to the top and bottom of the figure. This shows that a misspecification of \eq{\Psigmao^2\Plambok} --- whether it is overestimated or underestimated --- can only have a significant impact on the analysis error if \eq{\sigmao^2\lambok} and \eq{\sigmab^2\lambbk} are relatively close to each other at this spatial scale. For example, if \eq{\sigmao^2\lambok < \sigmab^2\lambbk}  by a factor of 10000, then overestimating \eq{\Psigmao^2\Plambok} by a factor as large as 100 has no significant impact on the analysis error (point indicated by the cross symbol in the figure). In such a case, the analysis will remain extremely close to the observations at this spatial scale. \par

%% file: experimental_framework.tex
The insight from the idealised problem in the previous section will be helpful for interpreting the results from the ocean data assimilation experiments presented in Section~\ref{sec_results}. In this section, we describe the ocean data assimilation system, focussing on the method for accounting for correlated error in along-track altimeter observations, which is a new development to the system. We also detail the experimental framework, which involves the assimilation of sea-surface height (SSH) innovations that are constructed from simulated error samples.

\subsection{The ocean data assimilation system}
\label{sec_model_and_da}

The experiments are performed with NEMOVAR \mbox{\citep{Mogensen_2009, Mogensen_2012}}, an incremental variational ocean data assimilation system that has been developed for the
ocean model NEMO \mbox{\citep{Madec_2023}}. NEMOVAR is used operationally at the European Centre for Medium-Range Weather Forecasts (ECMWF) and the UK Meteorological Office (Met Office; \textit{e.g.}, see \mbox{\cite{Chrust_2025}} and \mbox{\cite{Mignac_2024}} for up-to-date descriptions of the ECMWF and Met Office ocean data assimilation systems, respectively). For the purposes of this study where the innovation vector is simulated, only the numerical code for the quadratic cost function minimization (inner loop) of NEMOVAR is required.
The inner-loop configuration is global and similar to that of \cite{Balmaseda_2013}.
The global grid is the tri-polar (ORCA) grid of NEMO, with approximately $1^{\circ}$ horizontal resolution in the Extratropics and refined meridional resolution in the Tropics, with a minimum value of $0.3^{\circ}$ directly at the Equator. The model has $42$ vertical levels, 18 of which are in the upper 200~m.
Assimilation increments are computed for the ocean model state variables consisting of potential temperature (T), practical salinity (S), the two components of horizontal velocity (U and V) and SSH.

The data assimilation method is based on 3D-Var with the First-Guess at Appropriate Time (FGAT) approach. The control variables in the 3D-Var are T, and the so-called {\it unbalanced} parts of S and SSH. The control variables are linked to the state variables (T, S, U, V and SSH) via a non-linear balance operator \citep{Weaver_2005}. The background-error covariance matrix is taken to be block-diagonal with respect to the control variables. Each covariance block is factored into a symmetric product of a correlation matrix and a diagonal matrix of standard deviations.
The correlation matrices are modelled using a diffusion operator \citep{Weaver_2001}, with diffusion parameters chosen such that the  smoothing kernel (correlation function) is approximately Gaussian. The  diffusion-based operator approach for modelling correlation functions will be presented in more detail in the next section in the context of \eq{\bR}. For all three control variables, the background-error standard deviations and correlation length scales are parameterised in terms of either the background state or geographical position, as described in Appendix~B of \cite{Balmaseda_2013}. The specific values of the covariance parameters are given in Table~B1 of \cite{Balmaseda_2013}. %

Only SSH innovations are assimilated in the experiments since our main objective in this study is to evaluate the impact of the newly developed correlated-error model for sea-level anomaly (SLA) observations derived from along-track nadir altimeters. For generating the synthetic SSH innovations (see Section~\ref{sec:syn_da}), we only require the geographical locations of the SLA observations. For this, we use a realistic network of SLA observations from ENVISAT and Jason-2 altimeter observations on January 1st, 2009. Locations close to land and poleward of $50^{\circ}$S and $50^{\circ}$N have been rejected. The total number of SLA observation locations retained is 15460, which provides a dense coverage of the ocean surface. The SSH observation-error standard deviations are set to 5~cm and are inflated near coastlines as described in \cite{Mogensen_2012}.

Model estimates of SLA (\textit{i.e.}, the model-estimated SSH minus an estimate of the Mean Dynamic Topography) are mapped to the observation locations using bi-linear interpolation. The balance relations are linearised around a realistic background state from January 1st, 2009. The composition of the spatial interpolation operator and the linearised balance operator constitutes the linearised generalised observation operator (\eq{\bG}) that links the control variable increments to the innovation vector in Equation~\eqref{eq:Gapprox}, and that will be used to define the innovation itself in Equation~\eqref{eq_innov_error}. The quadratic cost function is minimised iteratively using the Restricted \eq{\bB}-preconditioned Conjugate Gradient (RBCG) algorithm described in \cite{Gurol_2014} (see their Algorithm~3). On each iteration of RBCG, methods are required to compute a matrix-vector product with \eq{\bB} and \eq{\bR^{-1}} (methods to apply \eq{\bB^{-1}} and \eq{\bR} are not required, however).

\subsection{Accounting for along-track error correlations in altimeter observations}
\label{sec_implementation}

In this section, we describe the method that has been implemented in NEMOVAR to account for along-track error correlations in nadir\footnote{As opposed to wide-swath altimeters such as SWOT.} altimeter data. The choice of nadir altimeter data is motivated by several factors. First, they are widely available and already assimilated in NEMOVAR; second, they are considered to be affected by errors correlated along track, which could be due to representativeness error associated with spatial scales observed but not resolved by the model \citep[\textit{e.g.}, ][]{Storto_2019}, or due to instrumental error associated with, for example, imperfect geophysical corrections or contamination between satellite footprints \citep[\textit{e.g.}, see Section~1.8 of][]{Stammer_2018};  third, their 1D nature vastly simplifies the implementation of a correlated error model in NEMOVAR, even though much of the code development is generic and applicable to other observation types as well; and fourth, the interpretation of the experimental results in relation to the theoretical 1D analysis is made much easier.\par

We assume that errors affecting observations from different nadir instruments are not cross-correlated. The matrix \eq{\bR} is thus block diagonal where each block (\eq{\bR_{\rm nadir}}) corresponds to a specific nadir instrument.
This matrix and its inverse can be written in the standard factored form
\begin{align}
\bR_{\rm nadir} & = \boldsymbol{\Sigma} \bC \boldsymbol {\Sigma}
\label{eq_Rnadir} \\
\mbox{and}
\hspace{0.5cm}
\bR_{\rm nadir}^{-1} & = \boldsymbol{\Sigma}^{-1} \bC^{-1} \boldsymbol {\Sigma}^{-1}
\end{align}
where \eq{\boldsymbol{\Sigma}} is a diagonal matrix of standard deviations and \eq{\bC} is a correlation matrix. To define the correlation matrix and its inverse, we use the implicit diffusion operator approach adapted to unstructured meshes described in \cite{Guillet_2019}. As with the implicit diffusion operator applied to structured meshes (\textit{e.g.}, for modelling \eq{\bB} in NEMOVAR as described in \cite{Weaver_2016}), the correlation matrix and its inverse have the general form
\begin{align}
   \bC & = \bGamma \bD \bGamma
   \label{eq_Cnadir} \\
   \mbox{and} \hspace{0.5cm}
   \bC^{-1} & = \bGamma^{-1} \bD^{-1} \bGamma^{-1}
\end{align}
where \eq{\bD} (\eq{\bD^{-1}}) is the matrix representing a diffusion (inverse diffusion) operator, and
\eq{\bGamma} is a diagonal matrix of normalisation factors that ensures that the diagonal elements of \eq{\bC} are approximately equal to 1. The main requirement for applying the diffusion-operator approach to \eq{\bR_{\rm nadir}} is to adopt an appropriate spatial discretization method for defining \eq{\bD} and \eq{\bD^{-1}} on a 1D mesh that can become unstructured along the satellite tracks due to, for example, missing observations or quality control decisions that result in observations being excluded from the assimilation.

\subsubsection{Finite-Element Method formulation}

Following \citet{Guillet_2019}, we discretize the diffusion operator using a finite-element method (FEM) with \eq{\mathcal{P}_1}-FEM basis functions (\textit{e.g.} see \mbox{\cite{Ciarlet_2002}}). We focus on the main points in this section and relegate the details on the numerical implementation to Appendix~\ref{app_implementation}. The symmetric expressions for \eq{\bD} and \eq{\bD^{-1}} with the FEM formulation are (cf. Equations~(42) and (43) in \citet{Guillet_2019})
\begin{align}
   \bD & = \left[ \left(\bM + \bA \right)^{-1} \bM \right]^{m} \bM^{-1}   \label{eq_D} \\
\mbox{and} \hspace{.5cm} \bD^{-1} & = \bM \left[ \bM^{-1} \left(\bM + \bA \right) \right]^{m},
\label{eq_inv_D}
\end{align}
where \eq{\bM} is the mass matrix, \eq{\bA} is the stiffness matrix\footnote{To denote the stiffness matrix, we use the notation \eq{\bA} instead of the standard notation \eq{\bK} (\textit{e.g.}, as presented in \mbox{\citet{Guillet_2019}}) to avoid confusion with the Kalman gain matrix.} and \eq{m} is the number of implicit diffusion iterations.
The mass and stiffness matrices are SPD and sparse \mbox{\citep[see Figure~4 of][]{Guillet_2019}} as the only non-zero off-diagonal elements are those associated with consecutive observation locations along the satellite track. In the case of nadir altimeter data, there are thus only three non-zero elements on each row of \eq{\bA}.  Matrix-vector products with \eq{\bM} and \eq{\bA} can thus be implemented efficiently even with a large volume of observations. Moreover, as observations only interact with their direct neighbours on the mesh, the matrix-vector products are straightforward to parallelize on multiple processors. In NEMOVAR, vectors in observation space can be distributed following the tiled model domain decomposition. Nodes located within each tile (processor) are considered local nodes. In addition to these local nodes, a halo is created that contains those nodes on other tiles which are connected to the local nodes by an edge of the mesh. The value of vectors inside this halo are updated once before each matrix-vector product with \eq{\bM + \bA} using Message Passing Interface (MPI) exchanges. \par

Let \eq{\Omega} denote the domain of the observations (in the case of nadir altimeters, the domain is the satellite track). Let the vector \eq{\bm{z}} denote a position on this domain. Let \eq{\psi_i} denote a \eq{\mathcal{P}_1}-FEM basis function of the \eq{i}-th observation, \textit{i.e.}; the function is equal to 1 at the position of the \eq{i}-th observation, equal to 0 at the position of other observations, and is linear in-between \citep[see Figure~3 of][]{Guillet_2019}. The \eq{kl}-th matrix element of \eq{\bM} and \eq{\bA} is given by the integral equations (cf. Equations~(27) and (28) in \mbox{\cite{Guillet_2019}})
\begin{align}
   \bM_{k,l} &= \int_{\Omega} \psi_k(\bm{z}) \psi_l(\bm{z}) d\bm{z} \label{eq_def_M} \\
\mbox{and} \hspace{.5cm}   \bA_{k,l} &= \int_{\Omega} \left[\bm{\kappa}(\bm{z}) \bm{\nabla} \psi_k(\bm{z})\right] ^{\rm T} \bm{\nabla}\psi_l(\bm{z}) d\bm{z}  \label{eq_def_A},
\end{align}
where \eq{\bm{\nabla}} denotes the gradient operator, and \eq{\bkappa} is the diffusion tensor. Equations~\eqref{eq_def_M} and \eqref{eq_def_A} show that each element of the matrices of \eq{\bM} and \eq{\bA}, respectively, can be expressed as an integral over the whole domain. Computing the value of these integrals separately for each matrix element is not computationally efficient. Instead, we use the partition of the domain into sub-domains by the mesh where the sub-domains correspond to the segments between observation locations in 1D (triangles in 2D). We cycle through each sub-domain, and increment all integrals that are non-zero on it before moving on to the next sub-domain. This process is usually referred to as the assembly process in the FEM literature (see Appendix~\ref{app_implementation} for more details). \par

To limit the cost of applying \eq{\bM^{-1}},
we adopt here the classical mass-lumping approach where \eq{\bM} is approximated by a diagonal matrix (see Appendix~\ref{app_implementation}). While this approximation increases the spatial discretization error, we did not observe any significant differences in the solution when compared to that obtained with the un-approximated mass matrix.
With mass lumping, applying \eq{\bD^{-1}} to a vector simplifies to \eq{m} products with a tri-diagonal matrix, and \eq{m+1} products with diagonal matrices. Applying \eq{\bR^{-1}} then requires 4 additional products with the diagonal matrices \eq{\boldsymbol{\Gamma}^{-1}} and \eq{\boldsymbol{\Sigma}^{-1}}.\par

For the specific case of 1D observations from nadir altimeter data, \eq{\bkappa} becomes a \eq{1\times 1} tensor, and thus a scalar diffusion coefficient. The diffusion coefficient controls the spatial extent of the along-track correlations. While it is straightforward to make the diffusion coefficient spatially varying, for the experiments in this article we used a constant diffusion coefficient denoted by the scalar \eq{\kappa}. The diffusion coefficient and the number of diffusion iterations (\eq{m}) are parameters that can be used to fit the correlation model to actual estimates of observation-error correlations. The effect of these parameters on determining the shape and spectral characteristics of the underlying correlation function is discussed in the next section. \par

\subsubsection{Specification of the diffusion-model parameters}

In the Euclidean space \eq{\R}, the smoothing kernel of the implicit diffusion operator with a constant diffusion coefficient is a covariance function from the Mat\'ern class \citep{Guttorp_2006, Weaver_2013}. Specifically, with \eq{m} implicit diffusion iterations, the covariance function is an Auto-Regressive (AR) function of order \eq{m}; \textit{i.e.}, a function characterised by a polynomial of order \eq{m-1} times the exponential function (\textit{e.g.}, see \cite{Mirouze_2010}). This analytical result provides the basis for controlling and interpreting the spatial and spectral smoothness properties of the implicit diffusion-based correlation model.\par

\begin{figure}[htb]
   \centering
  \includegraphics[width=\textwidth]{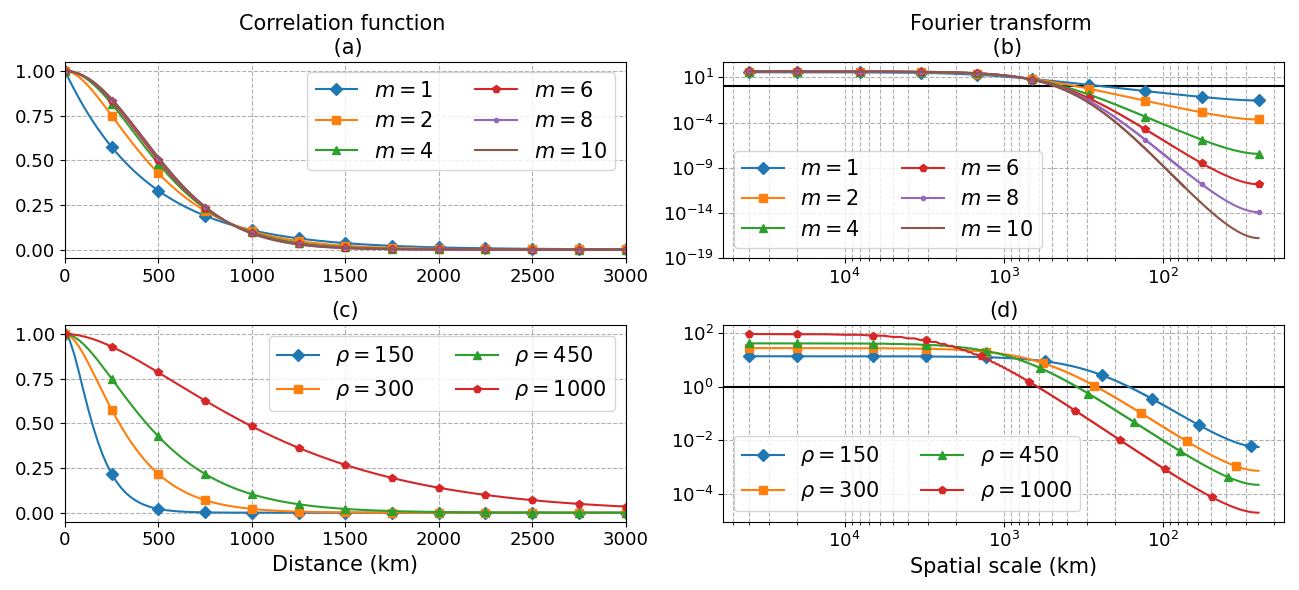}
   \caption{AR correlation functions (panels~(a) and (c)) and their respective Fourier transforms (panels~(b) and (d)). In panels~(a) and( b), \eq{\rho} is fixed to a value of 450~km and curves are displayed for different values of values of \eq{m}. In panels~(c) and (d), \eq{m} is fixed to a value of 2 and curves are displayed for different values of \eq{\rho}.}
   \label{fig_matern}
\end{figure}

The number of iterations \eq{m} controls the smoothness properties of the associated correlation function, as illustrated in Figure~\ref{fig_matern}. With \eq{m} large (values beyond about 10), the correlation function is approximately Gaussian (brown curve in Figure~\ref{fig_matern}(a)). Small values of \eq{m} correspond to correlation functions with a sharper (slower) decay at short (long) range. This is most clearly visible from the blue curve in Figure~\ref{fig_matern}(a), which corresponds to \eq{m=1}. For larger values of \eq{m}, the difference between the curves is difficult to distinguish in physical space, but the effect of \eq{m} is clearer in spectral space (Figure~\ref{fig_matern}(b)). In spectral space, \eq{m} is seen to control the steepness of the slope (or roll-off) of the spectrum at small spatial scales, with large values of \eq{m} (quasi-Gaussian functions) providing steepest roll-off. Background-error correlation functions are often specified with large values of \eq{m}; \textit{e.g.}, \eq{m} is set to 10 in the diffusion-based correlation models used for \eq{\bB} in the operational implementations of NEMOVAR at ECMWF and Met Office. One of the reasons motivating this choice is that the regularity of quasi-Gaussian functions allows them to be differentiated, which is typically required in the specification of multivariate cross-covariances \citep{Daley_1991,Derber_1999,Weaver_2005}. The constraint of regularity is less stringent for observation-error correlations. In fact, diagnosed observation-error correlations are often quite sharp near their peak \citep{Michel_2018, Waller_2016, Waller_2016b}, which is a feature that AR functions can capture with small values of \eq{m}.

In practice, the diffusion coefficient is specified indirectly through a more convenient length-scale parameter. Different definitions of correlation length scale have been proposed in the literature. A commonly used definition for background-error correlations is the Daley length scale \citep{Daley_1991, Pannekoucke_2008, Weaver_2013}, which is related to the local curvature of the correlation function at zero separation. However, the Daley length scale is not defined for non-differentiable correlation functions, which we do not want to exclude for modelling observation-error correlations. In particular, the exponential function, which corresponds to \eq{m=1} in the AR family, is not differentiable at zero separation. As in \mbox{\cite{Goux_2024}}, we use the following length-scale parameter from \cite{Lindgren_2011}:
\begin{equation}
   \rho = \sqrt{\kappa (2m - 1) },
   \label{eq_rho}
\end{equation}
which is very similar to the Daley length scale except for small values of \eq{m}  (see \cite{Goux_2024} for a discussion). A property of this length scale is that, at a distance \eq{\rho} from the AR function peak, the correlation is approximately 0.1 for all values of \eq{m}. The AR correlation functions displayed in Figure~\ref{fig_matern}(c) illustrate how increasing \eq{\rho} increases the spatial range of the correlations. In spectral space (Figure~\ref{fig_matern}(d)), \eq{\rho} affects the spatial scale at which the spectrum starts decaying at small spatial scales, but, unlike \eq{m}, does not change the decay rate itself.

Note that the correlation length scales specified in \mbox{\citet{Balmaseda_2013}} for background error are based on the Daley definition, not Equation~\eqref{eq_rho}. In Section~\ref{sec_results}, we will use Equation~\eqref{eq_rho} for both observation and background error. For \eq{m=10}, the length scale computed from Equation~\eqref{eq_rho} is about 5\% larger than the Daley length scale \mbox{\citep{Goux_2024}}. In Section~\ref{sec_results}, we will use subscripts o and b on \eq{\rho} and \eq{m} to distinguish the observation- and background-error correlation parameters, respectively.

\subsubsection{Normalisation factors}
\label{sec:norm}

The normalisation factors in \eq{\boldsymbol{\Gamma}} ensure that the specified standard deviations of the covariance matrix are not contaminated by amplitude errors in the representation of the correlation functions. The normalisation factors should be specified as the inverse of the square root of the diagonal elements of the diffusion-based covariance matrix \eq{\bD}. Therefore, to define the normalisation factors, we require a method to estimate the diagonal of \eq{\bD}. Different procedures have been proposed to estimate normalisation factors for diffusion-based formulations of \eq{\bB} (see \mbox{\cite{Weaver_2021}} for a review) and these procedures are generally applicable to diffusion-based formulations of \eq{\bR} as well. For the nadir altimeter network, the normalisation factors can be computed exactly using an efficient implementation of the brute-force procedure, as described below. If the observation network changes from one assimilation cycle to the next, then the normalisation factors have to be recomputed on each cycle.

The brute-force procedure requires a method to apply \eq{\bD} or a `square-root' factor of \eq{\bD} to a canonical vector (see Equations~(18) and (19) in \mbox{\cite{Weaver_2021}}.) Recall that only \eq{\bD^{-1}} is needed for the  minimization algorithm (RBCG in our case). In order to apply \eq{\bD}, we need a computationally efficient method to apply \eq{(\bM + \bA)^{-1}} in Equation~\eqref{eq_D}. Drawing on experience from the procedure used in \eq{\bB} in NEMOVAR \mbox{\citep{Weaver_2016}}, we use the Chebyshev Iteration (CI) to define an approximate inverse. Exploiting the fact that the mass matrix is diagonal with  mass lumping, which we denote by the matrix \eq{\bMl}, and hence can be easily factored into square-root  components \eq{\bMl^{1/2} \bMl^{1/2}}, we first reformulate \eq{\bD} in Equation~\eqref{eq_D} as
\begin{align}
   \bD & = \bMl^{- 1/2} \left[ \left(\bI + {\bMl}^{-1/2} \bA
   \bMl^{-1/2}\right)^{-1} \right]^{m}
   \bMl^{-1/2}.
   \label{eq_D_alt}
\end{align}

CI is used to solve the linear system involving the non-dimensional  matrix \eq{\bI + \bMl^{-1/2}\bA\bMl^{-1/2}}. From Equation~\eqref{eq_D_alt}, this system must be solved \eq{m} times, in sequence, where the solution of the linear system at a given iteration in the sequence is used as the right-hand side and first guess for the next linear system in the sequence. Since each system is solved only approximately, exact symmetry of the resulting (approximate) \eq{\bD} is not guaranteed. To avoid possible symmetry errors, we restrict \eq{m} to be an even number and factor Equation~\eqref{eq_D_alt} as
\begin{align}
   \bD & = \bMl^{-1/2} \left[\left(\bI + \bMl^{-1/2}\bA\bMl^{-1/2}\right)^{-1} \right]^{m/2} \left[\left(\bI + \bMl^{-1/2}\bA\bMl^{-1/2}\right)^{-1} \right]^{m/2}
   \bMl^{-1/2}.
   \label{eq_D_alt2}
\end{align}
Then, we use CI with a fixed number of iterations to solve the \eq{m/2} linear systems in one of the factors and the adjoint of the CI algorithm with the same number of iterations to solve the \eq{m/2} linear systems in the other factor.  This procedure preserves symmetry even in the presence of substantial truncation errors (see Section~5.2 in \cite{Weaver_2016} and Section~3.5 in \cite{Weaver_2018} for a more in-depth discussion). \par

CI requires as input an estimate of the extreme eigenvalues of the non-dimensional matrix \eq{\bI + \bMl^{-1/2}\bA\bMl^{-1/2}}. Following \citet{Weaver_2016}, we estimate the extreme eigenvalues prior to the RBCG minimization using a combined CG-Lanczos algorithm. This procedure does not add substantial cost to the initialization step in NEMOVAR.
After estimating the eigenvalues, the CI algorithm is used prior to the RBCG minimization to solve \eq{\bD \mathbf{z} = \mathbf{b}} with a random right-hand side in order to estimate the number of iterations required to reach a target tolerance of \eq{10^{-2}} for the residual 2-norm divided by the 2-norm of the right-hand side. Although this tolerance is larger than the tolerance of \eq{10^{-3}} suggested by \cite{Weaver_2016}, it is sufficient to achieve an average relative error of \eq{10^{-4}} for the normalisation factors. This procedure determines the number of iterations that are fixed when applying \eq{\bD} as explained above. In addition to ensuring symmetry, this procedure removes the need for estimating the norm of the residual on each iteration, and thus the need for global communications between different processors.

The brute-force approach requires the application of \eq{\bD} to a canonical vector that targets each of the observation locations. Although this method is generally not feasible for \eq{\bB} due to the large number of control variable grid points, it can be made efficient for the nadir altimeter network. Since observation points separated by a distance that is large compared to the correlation length scale do not interact significantly, we can estimate multiple normalisation factors simultaneously by applying \eq{\bD} to a vector with entries equal to 0 everywhere except at regularly spaced locations separated by (at least) \eq{r} times the correlation length scale where the entry is set to 1. In the case of nadir altimeters, empirical results indicate that a value of \eq{r=5} is sufficient to obtain an average relative error in the normalisation factors of \eq{10^{-4}}, which is an order of magnitude smaller than the error obtained with the standard method used for \eq{\bB}, which involves randomization with 10000 samples. With an observation-error correlation length scale of 125~km, \eq{r=5}, and altimeter observations separated by 20~km on average, this method requires only 32~applications of \eq{\bD} (compared to 10000 with randomization).

\subsection{Synthetic data assimilation experiments}
\label{sec:syn_da}

As described in this section, we use an experimental framework involving synthetic data to evaluate the impact of different representations of \eq{\bR}.
By combining Equations~\eqref{eq_def_yo} and \eqref{eq_def_innovation}, we can write the innovation vector as
\begin{equation}
   \bd = \G{\bxt} + \bepso - \G{\bxb},
   \label{eq_innov_decomp}
\end{equation}
where, from Equation~\eqref{eq_def_xb},
\begin{equation}
   \mathcal{G}(\bxt) = \mathcal{G}(\bxb - \bepsb) \simeq \mathcal{G}(\bxb) - \bG\bepsb,
   \label{eq_lin_Gxt}
\end{equation}
\eq{\bG} being the Jacobian matrix defined in Equation~\eqref{eq:Gapprox}.
Substituting Equation~\eqref{eq_lin_Gxt} in Equation \eqref{eq_innov_decomp} yields
\begin{equation}
   \bd \simeq \bepso - \bG\bepsb,
   \label{eq_innov_error}
\end{equation}
which shows that, in the linear approximation, the innovation vector can be expressed entirely in terms of the observation and background error.

In the synthetic data assimilation experiments, we use samples of background and observation error to simulate the innovation vector via Equation~\eqref{eq_innov_error}. To generate error samples with statistics consistent with \eq{\bB} and \eq{\bR}, we make use of `square-root' factors  \eq{\bU} and \eq{\bV} such that
\begin{equation}
   \bB = \bU \bU^{\rm T} \hspace{0.5cm}
   \mbox{and} \hspace{0.5cm}
   \bR = \bV \bV^{\rm T}.
\end{equation}
In NEMOVAR, both \eq{\bB} and \eq{\bR} are constructed directly in this way, using a product of operators. For example, the factor \eq{\bV} can be easily identified from Equations~\eqref{eq_Rnadir}, \eqref{eq_Cnadir} and \eqref{eq_D_alt2}.
Given \eq{\bU} and \eq{\bV}, we can generate an unlimited number of error samples, \eq{\ell = 1, 2, \ldots}, from the equations
\begin{equation}
   \bepso^{(\ell)} = \bV\bm{\eta}_{\rm o}^{(\ell)} \hspace{0.5cm} \mbox{and} \hspace{0.5cm}
   \bepsb^{(\ell)} = \bU\bm{\eta}_{\rm b}^{(\ell)},
\end{equation}
where \eq{\bm{\eta}_{\rm o}^{(\ell)}\in \mathbb{R}^p} and  \eq{\bm{\eta}_{\rm b}^{(\ell)}\in \mathbb{R}^n} are random vectors containing uncorrelated, Gaussian fields with mean equal to zero and standard deviation equal to one.

In all of our experiments, the background-error covariance matrix used to generate the error samples is the same as the one that is specified in the data assimilation algorithm. It has been described in Section~\ref{sec_model_and_da}.
However, for the observation-error covariance matrix,  we make a distinction between the `true' observation-error covariance matrix \eq{\bR}, which is used to generate the error samples, and the specified observation-error covariance matrix \eq{\PR}, which is used in the data assimilation algorithm (where we have used notation in alignment with that of Section~\ref{sec_theoretical}). Both \eq{\bR} and \eq{\PR} are defined according to the same observation-location mesh, but are specified with potentially different covariance parameters (\textit{i.e.}, \eq{\rhoo\neq\Prhoo}, \eq{\mo\neq\Pmo}, \eq{\sigmao\neq\Psigmao}).

As we use synthetic data, we can compute the actual analysis error when using different specifications of \eq{\bR} and \eq{\PR}. The analysis error can be expressed in terms of the analysis increment and background error as
\begin{equation}
   \bepsa \, = \, \bxa - \bxt
  \, = \, \bxa - \bxb + \bxb - \bxt
  \, = \, \bdxa + \bepsb.
  \label{eq_ana_err_min}
\end{equation}
As shown in \mbox{\citet{Goux_2024}}, different specifications of \eq{\bR} and \eq{\PR} can have a strong influence on the convergence of a $\mathbf{B}$-preconditioned CG minimization algorithm, like RBCG. Therefore, as a robust performance measure, we use Equation~\eqref{eq_ana_err_min} to track the evolution of the analysis error as a function of the iteration counter using the estimate of the increment on each iteration.  At each iteration of RBCG, the Root Mean Square (RMS) of the analysis error is computed for each state variable to assess both the convergence rate of the minimization and the accuracy of the solution. The convergence criterion for all experiments is defined by a 10-order of magnitude reduction in the \eq{\bB}-norm of the gradient relative to the \eq{\bB}-norm of the initial gradient.

%% file: small_corr.tex
\subsubsection{Expectations from theory}
\label{sec_small_corr_theory}

\begin{figure}[htb]
   \centering
   \includegraphics[width=\textwidth]{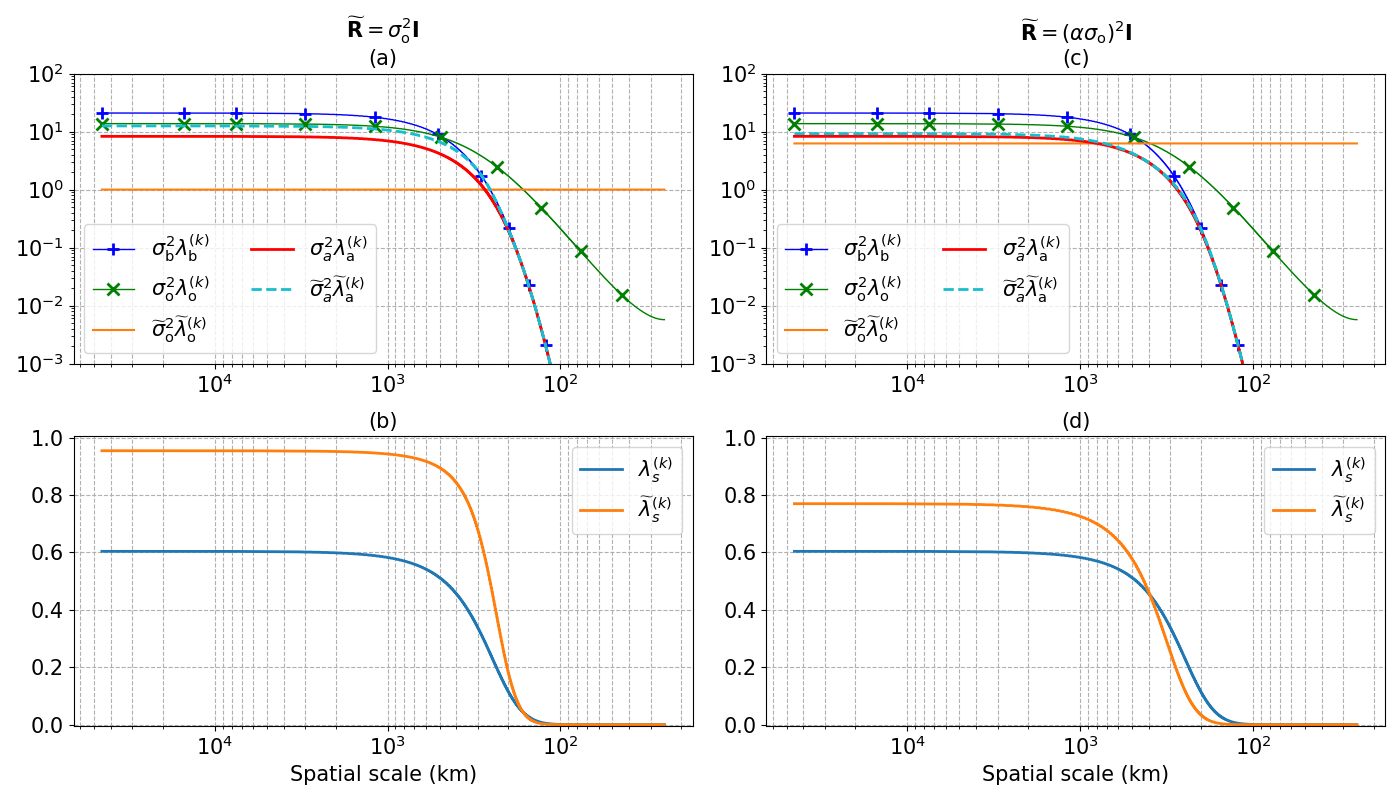}
   \caption{Panels~(a) and (c): scale-dependent variance of the background error (\eq{\sigmab^2\lambbk}),  observation error (\eq{\sigmao^2\lambok}) and its misspecified counterpart (\eq{\Psigmao^2\Plambok}),
   and analysis error (\eq{\sigmaa^2\lambak}; Equation~\eqref{eq_analysis_error_spectral}) and its sub-optimal counterpart (\eq{\Psigmaa^2\Plambak}; Equation~\eqref{eq_analysis_error_spectral_subopt}) as a function of spatial scale. Panels~(b) and (d): scale-dependent sensitivity (\eq{\lambsk}; Equation~\eqref{eq:lambda_sk}) and its sub-optimal counterpart (\eq{\Plambsk}; Equation~\eqref{eq_sensitivity_spectral_subopt}) as a function of spatial scale. The experiment is as in Figure~\ref{fig_scale_decomp_1}. The background-error correlation function is approximately Gaussian with a length scale \eq{\rhob=210}~km; the `true' observation-error correlation function is a SOAR function with a length scale \eq{\rhoo=150}~km; the misspecified  observation-error covariance function has no correlation, \textit{i.e.}, \eq{\PR=\sigmao^2\bI} in panels~(a) and (b), while \eq{\PR=\alpha^2\sigmao^2\bI}
   in panels~(c) and (d) where \eq{\alpha=2.5} is an optimally derived inflation factor that minimizes the analysis-error variance.}
   \label{fig_scale_decomp_L150}
\end{figure}

Figure~\ref{fig_scale_decomp_L150} summarizes the effect of neglecting observation-error correlations in the simplified framework presented in Section~\ref{sec_theoretical}. Here, the background-error correlation function is approximately Gaussian as in NEMOVAR (\eq{\mb=10}). The background-error correlation length scale is \eq{\rhob=210}~km. The (true) correlation length scale of the observation error is \eq{\rhoo=150}~km. The observation-error correlation function is a Second Order Auto Regressive (SOAR) function (\textit{i.e.}, \eq{\mo=2} in the implicit diffusion model). The scale-dependent sensitivity and analysis-error variance that would be obtained with optimal parameters (as in Figure~\ref{fig_scale_decomp_1}) are displayed as a reference. In addition, Figure~\ref{fig_scale_decomp_L150} displays the sub-optimal (practical) scale-dependent sensitivity (bottom panels) and analysis-error variance (top panels) obtained when using a diagonal \eq{\PR} in the gain matrix. In Figures~\ref{fig_scale_decomp_L150}(a) and (b), the correlations are simply ignored (\eq{\PR = \sigmao^2 \bI}). In this case, the observation-error variance is underestimated at large spatial scales and overestimated at small spatial scales (orange curve in Figure~\ref{fig_scale_decomp_L150}(a)). As a consequence, observations are overfit at large spatial scales and underfit at small spatial scales. Overfitting at large spatial scales translates into an excessively large sensitivity (left portion of Figure~\ref{fig_scale_decomp_L150}(b)) and a significant increase of the analysis-error variance at large spatial scales (\eq{\Psigmaa^2\Plambak>\sigmaa^2\lambak}; left portion of Figure~\ref{fig_scale_decomp_L150}(a)).\par

At small spatial scales, \eq{\Psigmao^2 \Plambok} is overestimated. The sub-optimal sensitivity thus becomes smaller than the optimal sensitivity at smaller spatial scales, but only at the smallest scales where both are close to zero (they overlap in the figure).
At the smallest spatial scales, \eq{\sigmab^2 \lambbk} is much smaller than both \eq{\Psigmao^2 \Plambok} and \eq{\sigmao^2 \lambok} (right portion of Figure~\ref{fig_scale_decomp_L150}(a)), and there is thus no significant increase of \eq{\Psigmaa^2 \Plambak} compared to the optimal case: in both the optimal and sub-optimal cases, the innovations (and hence the observations) are effectively discarded. \par

Overfitting observations at large spatial scales has a stronger impact than underfitting observations at small spatial scales on the analysis-error variance. This characteristic of the analysis justifies the common approach adopted in operational systems of inflating the observation-error variance; \textit{i.e.}, by setting \eq{\PR=\Psigmao^2\bI = (\alpha\sigmao)^2\bI} with \eq{\alpha > 1}. Although the choice of \eq{\alpha} is usually tuned empirically in operational configurations, in our idealized system, we can determine an optimal value of \eq{\alpha} that results, on average, in the minimum total analysis-error variance\footnote{Here, the average is computed from a sample of 100 analysis increments using innovation vectors generated randomly based on Equation~\eqref{eq_innov_error}.}. With the optimal \eq{\alpha} (2.5 for this example), Figure~\ref{fig_scale_decomp_L150}(c) and (d) show the best possible results that can be obtained with a diagonal \eq{\PR}. Increasing \eq{\alpha} increases \eq{\Psigmao \Plambok} at all spatial scales and thus decreases \eq{\Plambsk} at all spatial scales (Figure~\ref{fig_scale_decomp_L150}(d)). This approach can thus mitigate overfitting observations at large spatial scales and thus bring the sub-optimal analysis-error variance closer to the optimal one at those scales (left portion of Figure~\ref{fig_scale_decomp_L150}(c)), but at the expense of worsening the fit at the small spatial scales. The optimal \eq{\alpha} is not large enough to increase significantly the impact of underfitting at small spatial scales on the analysis error (right portion of Figure~\ref{fig_scale_decomp_L150}(c)). \par

In summary, if \eq{\PR} is diagonal, variance inflation is advantageous for the case \eq{\rhoo < \rhob} as it greatly reduces the analysis-error variance at large spatial scales by reducing the fit to the observations. Importantly, it does so without substantially increasing the analysis-error variance at small spatial scales. Although inflation further reduces the fit to the observations at the small spatial scales, this has little impact on the analysis error. This small impact is explained by the fact that, with or without an inflated \eq{\sigmao}, the sensitivity is negligible at these small spatial scales since \eq{\sigmao^2 \lambok} is larger than \eq{\sigmab^2 \lambbk} by several orders of magnitude. However, variance inflation is less effective than actually accounting for observation-error correlations. With the optimal level of variance inflation, the sub-optimal sensitivity at large spatial scales remains too large compared to the optimal sensitivity. Increasing the inflation factor \eq{\alpha} further would reduce the overfitting at large spatial scales; however, in this case, the observations would be underfit at intermediate spatial scales as well as small spatial scales. Underfitting the observations at intermediate spatial scales would significantly degrade the analysis since, contrary to the small spatial scales, \eq{\sigmab^2 \lambbk} and \eq{\sigmao^2 \lambok} are comparable in magnitude at these scales. Although variance inflation can mitigate the detrimental effects from neglecting correlated error, accounting for observation-error correlations is still beneficial. This is assessed experimentally in NEMOVAR in the next section.

\subsubsection{NEMOVAR experiments }
\label{sec_small_corr_practical}

\begin{figure}[h!]
\begin{center}
   \includegraphics[width=0.5\textwidth]{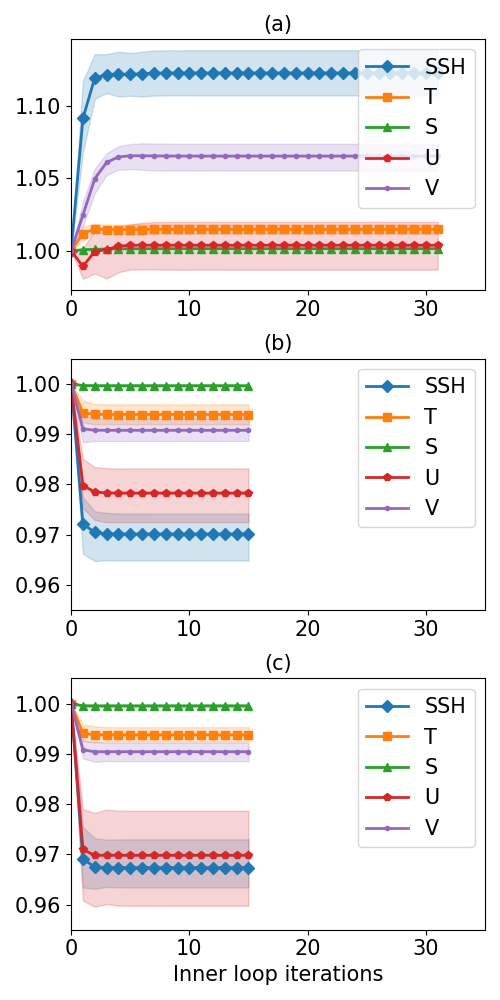}
   \caption{RMS of the analysis error at each iteration of the inner loop, for each variable (SSH: sea surface height, T: temperature, S: salinity, U: zonal current velocity, V: meridional current velocity). The values have been normalised by the RMS of the background error at each iteration. The background-error correlation function is approximately Gaussian (\eq{\mb=10}) where the length scale \eq{\rhob} varies between 210~km and 420~km. The true observation-error correlation function is a SOAR function (\eq{\mo=2}) with a length scale \eq{\rhoo =150}~km. Solid lines correspond to  the mean of 100 samples of the RMS of the analysis error derived from different realisations of the random innovation vector (Equation~\eqref{eq_innov_error}). Filled-in areas indicate the spread between the 25\% and 75\% quantiles of the samples. (a) \eq{\PR} is specified as a diagonal matrix; (b) \eq{\PR} is specified as a diagonal matrix with optimal inflation; and (c) \eq{\PR = \bR}, the true observation-error covariance matrix. For panel~(b), the optimal inflation factor is \eq{\alpha = 2.5}.}
   \label{fig_convergence_L150}
   \vspace{-0.5cm}
\end{center}
\end{figure}

We consider three experiments to determine if the conclusions deduced from the simplified system in Section~\ref{sec_small_corr_theory} carry over to a realistic configuration. Correlated SSH observation errors are simulated with the implicit diffusion model with parameter settings \eq{\rhoo=150}~km and \eq{\mo=2} to emulate the SOAR function of Section~\ref{sec_small_corr_theory}. For comparison, the smallest specified background-error correlation length scale in the meridional direction (which is approximately aligned with the satellite track) is 210~km \mbox{\citep{Balmaseda_2013}}. The simulated innovations are then assimilated with three different specifications of \eq{\PR}: a naive diagonal \eq{\PR} without variance inflation; an optimally inflated diagonal \eq{\PR}; and an accurate non-diagonal \eq{\PR=\bR}. The RMS of the analysis error for each variable at each iteration (Equation~\eqref{eq_ana_err_min}) is shown for the three \eq{\PR} specifications in panels (a), (b), (c) of Figure~\ref{fig_convergence_L150}, respectively. The RMS values are computed over all model grid points. The values have been normalized by the corresponding RMS of the background error. As we observed a strong variability in the RMS of the analysis error with different random realisations of the innovation vector, we compute the average RMS from a sample of analysis errors from 100 experiments. The solid lines in Figure~\ref{fig_convergence_L150} indicate the mean of the sample. The area between the 25\% and 75\% quantiles of each sample is filled in to give a qualitative indication of the spread of the sample.

If observation-error correlations are not accounted for (\eq{\PR=\sigmao^2\bI}), we expect a significant degradation of the analysis error due to the overfit of the observations at large spatial scales. This degradation appears very clearly in Figure~\ref{fig_convergence_L150}(a), where the analysis error is larger than the background error for all variables. The analysis error is particularly large for SSH. Figure~\ref{fig_convergence_L150}(b) shows that, as predicted in Section~\ref{sec_small_corr_theory}, variance inflation can successfully mitigate the overfitting of the observations at large spatial scales and results in an improved analysis for all variables. The minimization convergence is also accelerated by variance inflation, requiring 15~iterations compared to 31~iterations without inflation. As expected, increasing \eq{\Psigmao} decreases the condition number of the (\eq{\bB}-preconditioned) Hessian matrix of the cost function \citep{Tabeart_2021,Goux_2024}, in this case from 9.5 to 2.4. The relatively low values of the condition number are associated with the coarseness of the grid in this configuration. The sensitivity of each variable to variance inflation is discussed in more detail in Section~\ref{sec_inflation}.\par

 Finally, accounting for observation-error correlations with an accurate non-diagonal \eq{\PR=\bR} (Figure~\ref{fig_convergence_L150}(c)) leads to an even lower analysis error than with optimal inflation. Furthermore, this improvement occurs without a substantial degradation of the convergence rate, as the condition number is similar with optimal variance inflation and with a non-diagonal \eq{\PR} (2.4 and 2.5, respectively). This is consistent with theory where we expect the condition number of the Hessian matrix to remain similar or to be reduced when observation-error correlation length scales are short compared to background-error correlation length scales \citep{Goux_2024}. The improved analysis is due to the absence of overfitting the observations at large spatial scales and underfitting the observations at small spatial scales, which occurs when using a diagonal \eq{\PR}. The improvement at small spatial scales is highlighted by comparing the response of the different variables to the use of a non-diagonal \eq{\PR}. \par

The variable most improved by the assimilation of altimeter data is SSH as it corresponds directly to the measured quantity. The other variables also see an improvement due to the balance relations. Temperature and salinity are related to SSH via a T-S relation, linearised equation of state and dynamic height relation \citep{Weaver_2005}, although the effect on salinity is much smaller due to the smaller contribution of salinity  to the expansion/compression of the water column. At the surface, the analysis increment for the zonal velocity component (\eq{\delta U}) and meridional velocity component (\eq{\delta V}) are related to the horizontal gradient of the SSH increment (\eq{\delta \eta}) through the geostrophic balance relation \citep{Weaver_2005}:
\begin{align}
   \delta U &= -\frac{g}{f}\frac{1}{a}\frac{\partial \delta \eta} {\partial \phi}, \label{eq_geostrophy}\\
   \delta V &= \frac{g}{f}\frac{1}{a \cos \phi}\frac{\partial \delta \eta}{\partial \lambda}
\end{align}
where \eq{g} denotes standard gravity, \eq{f} is the Coriolis parameter, \eq{a} is the radius of the Earth, \eq{\phi} is latitude and \eq{\lambda} is longitude. The velocity components are thus more influenced by the \textit{gradient} of SSH than SSH itself. We can then expect them to improve especially when the small spatial scales of the SSH analysis are improved, which is the case when accounting for observation-error correlations \citep{Seaman_1977}. The improvement in the velocity components is indeed striking when comparing panels~(b) and (c) of Figure~\ref{fig_convergence_L150}. Although the zonal velocity is not the observed variable, the reduction in its RMS analysis error when using a non-diagonal \eq{\PR} instead of optimal inflation is more significant than that of SSH. Since the nadir satellite tracks are almost aligned with the meridional axis, the SSH meridional derivative (\eq{\partial \delta \eta / \partial \phi}) is better sampled than the zonal derivative (\eq{\partial \delta \eta / \partial \lambda}). Furthermore, as the observation errors are correlated along track, and that these are correctly accounted for, information about the short spatial scales (gradient) is more efficiently extracted from the observations in the meridional direction. This improvement in the retrieval of the meridional derivative of SSH translates into an improvement in the retrieval of zonal velocity as their respective increments are proportional via the geostrophic relation (Equation~\eqref{eq_geostrophy}). As the zonal derivative of SSH is almost orthogonal to the satellite track, these observations provide less information on the meridional component of the geostrophic velocity. \par

In summary, if observation-error correlations have short length scales compared to those of background-error correlations (here, \eq{\rhoo=150}~km, while \eq{\rhob} ranges from 210~km to 420~km in the meridional direction), not accounting for them can lead to a negative impact of the observations on the analysis due to, in particular, overfitting the observations at large spatial scales. In this situation, variance inflation can successfully reduce the overfitting at large spatial scales while simultaneously improving the convergence of the minimization. However, improvement at large spatial scales occurs at the expense of neglecting gradient information in the observations to improve the small spatial scales. On the other hand, using a non-diagonal \eq{\PR} allows the small-scale information to be extracted from the observations  without overfitting the observations at large spatial scales. Furthermore, this improvement can be achieved without degrading the convergence rate of the minimization.

%% file: large_corr.tex
In the previous section, we showed that using a non-diagonal \eq{\PR} has a positive impact on the analysis when \eq{\rhoo < \rhob}. In this section, we consider the opposite case when \eq{\rhoo > \rhob}. The background-error correlations are specified as in the previous section (\eq{\rhob} ranges from 210~km to 420~km in the meridional direction and \eq{\mb=10}), but now we set the observation-error correlation length scale to \eq{\rhoo=450}~km, while keeping the same value of the smoothness parameter (\eq{\mo = 2}).\par

\subsubsection{Expectations from theory}
\label{sec_large_corr_theory}

\begin{figure}[h]
   \includegraphics[width=\textwidth]{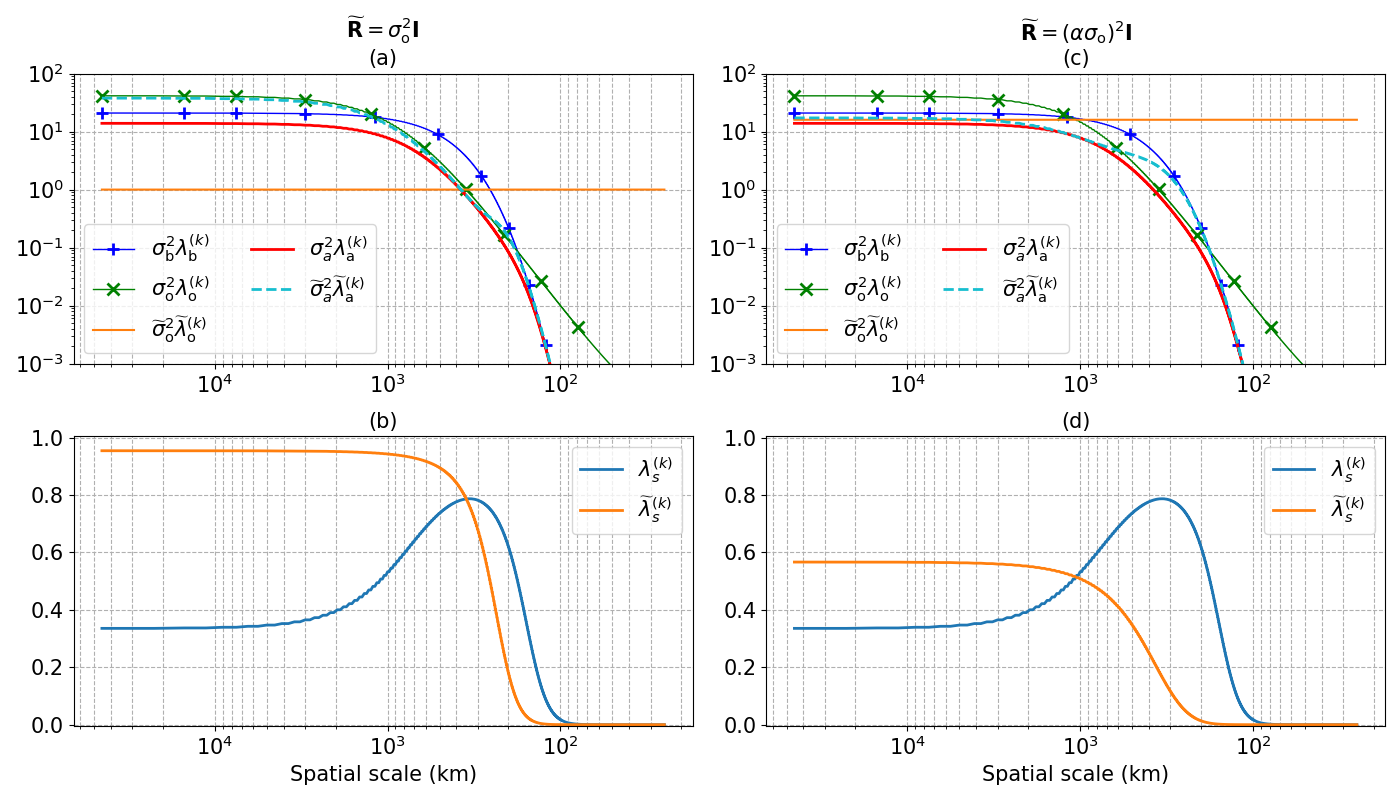}
   \caption{Same as for Figure~\ref{fig_scale_decomp_L150} but with \eq{\rhoo = 450}~km. The optimal inflation factor in this case is \eq{\alpha = 4}.}
   \label{fig_scale_decomp_L450}
\end{figure}

In Section~\ref{sec_small_corr}, the background was less accurate than the observations at large spatial scales and more accurate at small spatial scales. With a larger \eq{\rhoo}, this is no longer the case as illustrated in  Figure~\ref{fig_scale_decomp_L450}. At large spatial scales, \eq{\sigmao^2 \lambok} is slightly larger than \eq{\sigmab^2 \lambbk} (left portion of panels~(b) and (d)) and \eq{\lambsk} is smaller than 0.4 (left portion of panels~(a) and (c)), which leads to an analysis closer to the background than to the observations. However, since \eq{\rhoo > \rhob}, \eq{\sigmao^2 \lambok} begins to decay sooner than \eq{\sigmab^2 \lambbk} when approaching smaller spatial scales. This leads to an increase in \eq{\lambsk} at intermediate to small spatial scales (between approximately 2000~km and 200~km in this example) where observations become predominant in the optimal analysis.\par

With smaller \eq{\rhoo}, the shape of the optimal and suboptimal sensitivities \eq{\lambsk} and \eq{\Plambsk} were relatively similar as shown in Figure~\ref{fig_scale_decomp_L150}. This is not the case with larger \eq{\rhoo} as can be seen in Figure~\ref{fig_scale_decomp_L450}; with a diagonal \eq{\PR}, \eq{\Plambsk} always decreases monotonically when approaching smaller spatial scales, while \eq{\lambsk} (with a non-diagonal \eq{\PR}) shows a more complex behaviour. Consequently, the optimal sensitivity cannot be approximated correctly with a diagonal \eq{\PR}. In particular, there is no level of variance inflation that allows the analysis to extract relevant information from intermediate to small scales from the observations without significantly overfitting the observations at large spatial scales. This ineffectiveness of variance inflation in the presence of observation errors with large length scales has been observed in practice by \mbox{\cite{Reid_2020}} in the context of satellite-retrieved SST data assimilation.

\subsubsection{NEMOVAR experiments}
\label{sec_large_corr_practical}

As in Section \ref{sec_small_corr_practical}, we consider three sets of experiments involving different specifications of \eq{\PR}: a naive diagonal \eq{\PR} without variance inflation, an optimally inflated diagonal \eq{\PR}, and an accurate non-diagonal \eq{\PR=\bR}. The correlated SSH observation errors are simulated with the implicit diffusion model with parameter settings \eq{\rhoo=450}~km and \eq{\mo=2}. \par

\begin{figure}[h!]
\begin{center}   \includegraphics[width=0.5\textwidth]{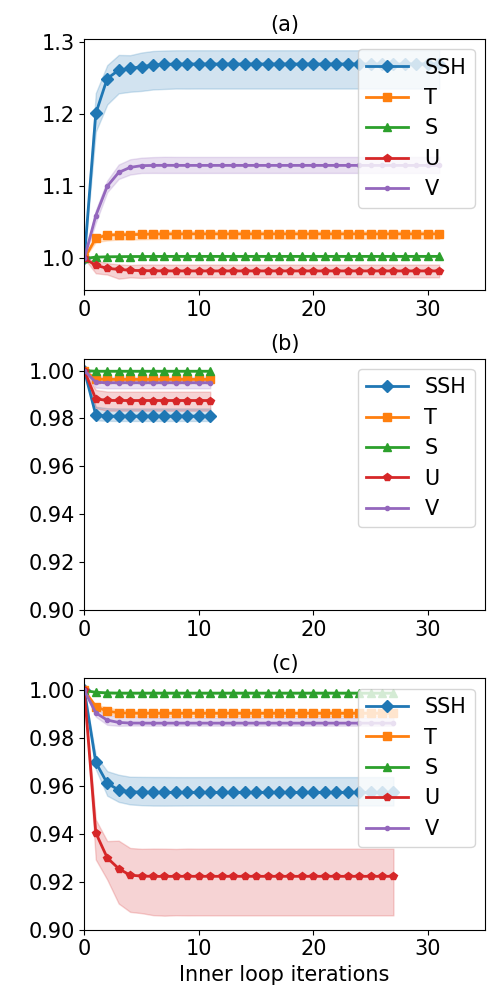}
   \caption{Same as for Figure~\ref{fig_convergence_L150} but with \eq{\rhoo = 450}~km. For panel~(b), the optimal inflation factor is \eq{\alpha = 4}.}
   \label{fig_convergence_L450}
\end{center}
\end{figure}

When the observation-error correlations are neglected, the analysis is degraded compared to the background, for all variables except zonal velocity, which achieves a 5\% reduction in RMS error (Figure~\ref{fig_convergence_L450}(a)). The difference between zonal velocity and the other variables can be explained by referring to panels~(a) and (b) of Figure~\ref{fig_scale_decomp_L450}. The analysis overfits the observations at large spatial scales, which can explain the degradation for most variables. The scale-dependent analysis-error variance appears to be particularly sensitive to overfitting at large spatial scales. This is consistent with the theory developed in Section~\ref{sec_misspec}, which showed that overfitting the least accurate source of information (here the observations) has more impact on the analysis error than underfitting it.

However, as discussed in Section~\ref{sec_small_corr_practical}, the zonal  velocity is influenced more by the meridional  derivative of the assimilated SSH observations than the absolute value of SSH, implying that it is mainly affected by the small spatial scales in the observations. We can see in Figure~\ref{fig_scale_decomp_L450}(b) that without variance inflation, \eq{\Plambsk} is relatively close to \eq{\lambsk} at small spatial scales, and in Figure~\ref{fig_scale_decomp_L450}(a) that the increase of analysis error at these scales is small (\eq{\sigmaa^2\lambak\simeq \Psigmaa \Plambak}). Consequently, we expect \eq{\PR=\sigmao^2\bI} to be inappropriate for most variables except for quantities such as velocity that depend on the horizontal gradient (or higher-order derivatives) of the observations. The same effect is not observed for meridional velocity, as it is related to the zonal derivative of SSH, which is almost perpendicular to the satellite track and thus less constrained by the observations.\par

With optimal variance inflation, overfitting at large spatial scales can be mitigated but not without preventing the extraction of small-scale information from the observations. Optimal inflation is particularly ineffective in this configuration as the small spatial scales are where the observations provide most of the information for improving the background (\textit{i.e.}, where \eq{\sigmao^2\lambok < \sigmab^2\lambbk}). Although optimal variance inflation can prevent a negative impact of the observations on the analysis, it also negates most of their positive impact (Figure~\ref{fig_convergence_L450}(b)). \citet{King_2021} reached a similar conclusion when assimilating simulated wide-swath altimeter observations with correlated error in a regional configuration. In their case, they had to apply a large amount of data thinning to compensate for correlated error when using a diagonal \eq{\PR}. While this prevented the analysis from being degraded by the simulated observations, it also significantly reduced their impact on the analysis compared to a case where the simulated observations were taken to have uncorrelated errors.

The impact of the observations is improved by accounting for the correlated error. With an accurate non-diagonal \eq{\PR = \bR}, the error reduction of the SSH reaches 4\% compared to 2\% with optimal variance inflation (cf. Figure~\ref{fig_convergence_L450}(c) and (b)). The effect on zonal velocity is even more impressive, with the relative error reduction reaching about 8\%. It is interesting to note that, comparing the relative reduction of SSH and zonal velocity, the gradient of the observations provides more useful information than the observations themselves; \textit{i.e.},
SSH observations with correlated error such that \eq{\rhoo > \rhob} provide a stronger constraint on the geostrophic component of the zonal velocity than on the SSH itself. This result supports the approach  to assimilate SSH observations together with their derivatives \citep{Brankart_2009,Ruggiero_2016}.\par

Using the non-diagonal \eq{\PR} has a slight positive impact on the convergence rate of the minimization compared to the experiments with the naive diagonal \eq{\PR}, which uses the same observation-error variance. Although there is a risk of increasing the condition number with a non-diagonal \eq{\bR} when \eq{\rhoo > \rhob}, it is limited by the fact \eq{\mo<\mb} (see Figures~7 and 8 in \cite{Goux_2024}). In the case considered here, \eq{\rhoo} (equal to \eq{450}~km) is not large enough compared to \eq{\rhob} (which ranges from 210~km to 420~km) for the condition number to be degraded by the non-diagonal \eq{\PR}. The condition number is 9.2 with a non-diagonal \eq{\PR} compared to 9.5 with the non-inflated diagonal \eq{\PR}. An increase in the condition number and a degraded convergence rate do occur if the same experiment is performed with \eq{\rhoo=800}~km (not shown). In this case the condition number is 45 and full convergence is reached in 53~iterations.
On the other hand, the convergence with an optimally inflated \eq{\PR} is faster than with a non-diagonal \eq{\PR}, as full convergence is reached in 11~iterations compared to 28~iterations. The improved convergence rate is explained by the larger observation-error variance in the inflated \eq{\PR} experiment, which decreases the condition number to 1.5. Even though reaching full convergence requires more iterations with an inflated non-diagonal \eq{\PR}, most of the error reduction occurs in the early iterations. In fact, at each iteration, the RMS of the analysis error is lower than what can be achieved with any diagonal approximation.

%% file: app_inflation.tex
Optimal variance inflation allows the practical observation-error variance to achieve a better match to the true variance at large spatial scales and thus to reduce overfitting of the observations at those scales (Section~\ref{sec_theoretical}). However, this can only occur at the expense of the small spatial scales where the discrepancy between practical and true observation-error variance becomes worse. In the NEMOVAR experiments, zonal velocity was shown to be more sensitive than SSH to small spatial scales. In this section, we compare the analysis errors for these two variables to evaluate their response to different levels of variance inflation.\par

\begin{figure}[h!]
   \centering
   \includegraphics[width=\textwidth]{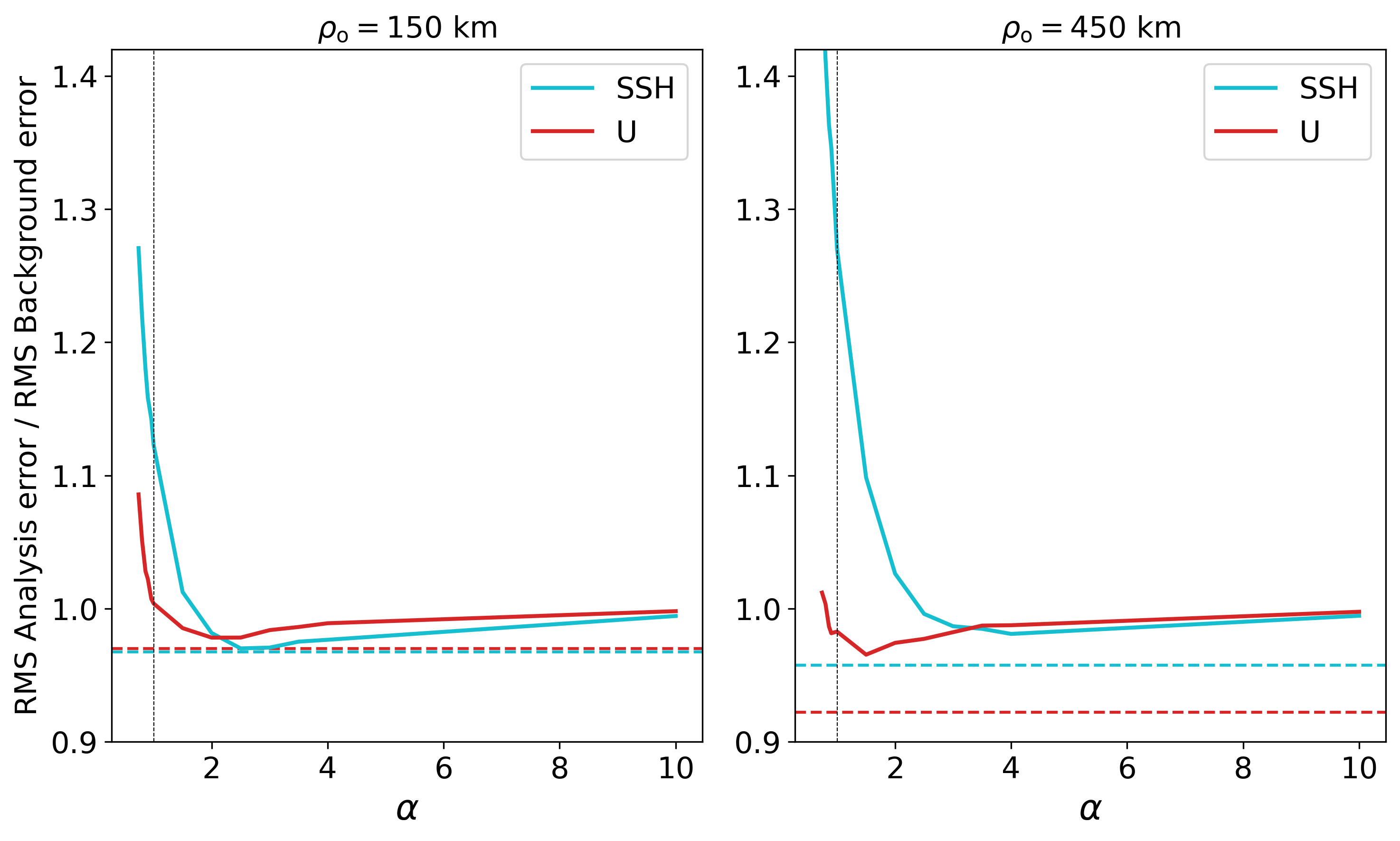}
   \caption{The RMS of the analysis error for SSH and zonal velocity for different values of the variance inflation factor \eq{\alpha}.  The results are averaged over 100 NEMOVAR experiments and have been normalised by the RMS of the background error. The true observation-error correlation function is a SOAR function with a length scale of 150~km in panel (a) and 450~km in panel (b). The practical \eq{\PR} is specified as \eq{(\alpha\sigmao)^2\bI}. The dotted horizontal lines indicate the values reached with a non-diagonal \eq{\PR=\bR}. The black vertical dashed line marks no inflation (\eq{\alpha  =1}).}
   \label{fig_inflation}
\end{figure}

If \eq{\rhoo=150}~km (\textit{i.e.}, smaller than \eq{\rhob} everywhere in the domain), variance inflation leads to an improvement for both variables (Figure \ref{fig_inflation}a). Compared to SSH, the zonal velocity seems to perform generally better with lower values of \eq{\alpha}. This is consistent with our previous conclusion that the small spatial-scale information in the observations (\textit{i.e.}, the meridional derivative of SSH and thus the  zonal velocity through geostrophy) is best extracted with less variance inflation. For zonal velocity, there is no level of inflation with a diagonal \eq{\PR} for which the RMS of the analysis error is as low as that obtained with a non-diagonal \eq{\PR} (marked by dashed lines in Figure~\ref{fig_inflation}). \par

This is even more striking when \eq{\rhoo=450}~km, as seen in Figure~\ref{fig_inflation}(b). Even if \eq{\alpha} is set optimally for a specific variable, the analysis error of this variable would still be much larger than it would be with a non-diagonal \eq{\PR}. Moreover, there is a large discrepancy between optimal inflation factors for SSH and zonal velocity, with \eq{\alpha =4} and \eq{\alpha=1.5} being the optimal values for these variables, respectively. Consequently, variance inflation forces a compromise on the accuracy of some variables in favour of others. There is no level of variance inflation that prevents the SSH analysis from being worse than the background without negating the improvement of the zonal velocity analysis.\par

%% file: stability.tex
The results from the NEMOVAR experiments illustrated the sensitivity of the analysis error and its convergence behaviour to the correlation length scale in \eq{\PR}. Additional sensitivity experiments were conducted to evaluate the impact from adopting a diagonal \eq{\PR} approximation when the actual observation errors are correlated. In this section, we examine the sensitivity from a different perspective, focussing on the non-diagonal \eq{\PR^{-1}} and its sensitivity to noise at small spatial scales.

First, consider the linear system \eq{\PR \mathbf{z} = \mathbf{b}} with \eq{\mathbf{z}, \mathbf{b}\in \mathbb{R}^{p}}, \eq{p} being the number of observations.
Let \eq{\delta\mathbf{z}} denote the perturbation to the solution \eq{\mathbf{z}} that results from a perturbation \eq{\delta\mathbf{b}} to the right-hand side \eq{\mathbf{b}}; \textit{i.e.}, \eq{\PR ( \mathbf{z} + \delta\mathbf{z}) = \mathbf{b} + \delta\mathbf{b}} where \eq{\delta \mathbf{z} = \PR^{-1} \delta \mathbf{b}}. The condition number of \eq{\PR}, denoted \eq{\kappa(\PR)}, bounds the ratio of the norm of \eq{\delta \mathbf{z}} to the norm of \eq{\mathbf{z}} \citep[Section~2.6.2 of][]{Golub_2013}:
\begin{equation}
   \frac{\|\delta\mathbf{z}\|_2}{\|\mathbf{z}\|_2} \leq \kappa(\PR) \frac{\|\delta\mathbf{b}\|_2}{\|\mathbf{b}\|_2}.
   \label{eq_cond_stab}
\end{equation}
A small perturbation to \eq{\mathbf{b}} will not result in a large perturbation to \eq{\mathbf{z}} if \eq{\kappa(\PR)} is small but could lead to a large perturbation to \eq{\mathbf{z}} if \eq{\kappa(\PR)} is large. Furthermore, in the latter case, the perturbation can be amplified by successive applications of \eq{\PR^{-1}} in the minimization algorithm. For example, if the experiment of Figure~\ref{fig_convergence_L450} is repeated with \eq{\Pmo = \mo = 10} instead of \eq{\Pmo = \mo = 2}, the noise inherent in finite-precision calculations can be amplified to an extent that the analysis becomes worse than the background even though \eq{\PR =\bR}. The degraded performance is explained by the huge condition number of \eq{\PR}, which is \eq{10^{9}} in this case.\par

A large condition number does not necessarily imply that perturbations will be excessively amplified (Equation~\eqref{eq_cond_stab} is only an upper bound), but the shape of the spectrum of \eq{\PR^{-1}} indicates how this extreme sensitivity can occur. The spectrum of \eq{\PR^{-1}} is the inverse of the spectrum of \eq{\PR}. When \eq{\PR^{-1}} is applied to a vector, it will amplify the small spatial scales of that vector by several orders of magnitude, with the smallest spatial scales amplified the most. In particular, referring to Figure~\ref{fig_matern}(b), we see that, for \eq{\mo=10} and \eq{\rhoo=450}~km, a condition number of \eq{10^{9}} would result in \eq{\PR^{-1}} damping the largest spatial scale by a factor of 10 while amplifying the smallest spatial scale by a factor of \eq{10^{8}}. This property of Gaussian-like functions makes them and other correlation functions with a sharp decay rate at small spatial scales inappropriate for modelling \eq{\bR^{-1}}. As expected from Equation~(17) in \mbox{\cite{Goux_2024}}, by decreasing \eq{\mo} from 10 to 2, the largest eigenvalue of \eq{\PR^{-1}} decreases while its smallest eigenvalue is virtually unchanged, resulting in a lower condition number (541 in the NEMOVAR experiment). With \eq{\mo = 2}, the small spatial scales are no longer excessively amplified, as confirmed by the positive results presented earlier in Figure~\ref{fig_convergence_L450}(c).

\begin{figure}[h!]
\begin{center}   \includegraphics[width=0.5\textwidth]{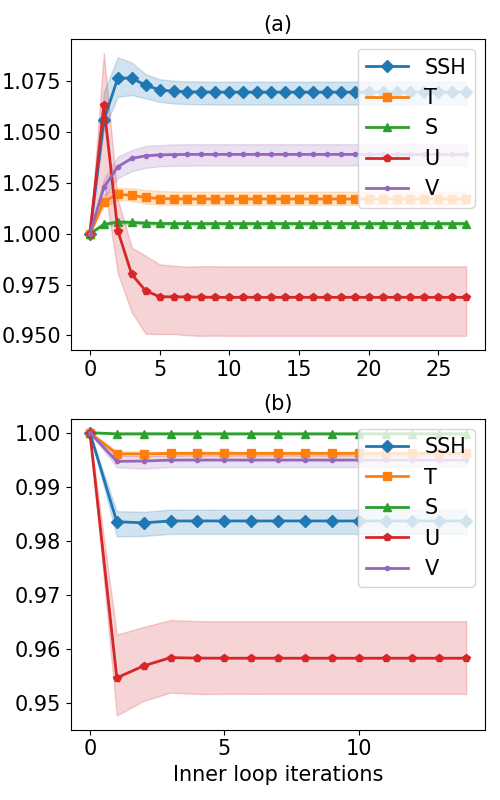}
   \caption{Same as Figure~\ref{fig_convergence_L450} but with the true (correlated) observation error including an additional component of uncorrelated error with variance equal to 1\% of the total observation-error variance (\eq{\sigmao^2 = 1}). (a) \eq{\rhoo = \Prhoo = 450}~km; (b) \eq{\rhoo = 450}~km, \eq{\Prhoo = 400}~km and \eq{\Psigmao = 2} (\textit{i.e.}, the inflation factor is \eq{\alpha = 2}).}
   \label{fig_convergence_unstab}
\end{center}
\end{figure}

In a real-data application, the observation-error covariances are not perfectly known, which can lead to an additional risk. Let us assume that the true observation error is the sum of a correlated component \eq{\bm{\epsilon}_{\rm corr}} and an uncorrelated component \eq{\bm{\epsilon}_{\rm uncorr}}:
\begin{equation}
    \bepso = \frac{1}{\sqrt{1+\beta^2}} \left( \bm{\epsilon}_{\rm corr} + \beta \bm{\epsilon}_{\rm uncorr} \right)
    \label{eq_composite_error}
\end{equation}
where \eq{\beta} is a scalar set to 0.01 so that the uncorrelated error represents 1\% of the total variance of the observation error. The normalisation factor \eq{\sqrt{1+\beta^2}} ensures that the total variance is equal to the variance of \eq{\bm{\epsilon}_{\rm corr}}. The correlation model parameters are set  as before: \eq{\rhoo=450}~km and \eq{\mo=2}. The specified \eq{\PR} ignores this additional uncorrelated component: \textit{i.e.}, \eq{\PR= \E {\bm{\epsilon}_{\rm corr}  \bm{\epsilon}_{\rm corr}^{\rm T}}} while \eq{\bR= \E {\bepso \bepso^{\rm T}}}. Even though the perturbation to the observation error is small compared to the total observation error, it leads to a significantly degraded analysis  where all fields except for zonal velocity are worse than the background (Figure~\ref{fig_convergence_unstab}(a)). This negative impact sharply contrasts the positive impact seen in Figure~\ref{fig_convergence_L450}(c) in the absence of this additional noise. The degradation comes from \eq{\PR^{-1}} excessively amplifying the errors at small spatial scales where the actual observation errors contain slightly more power than accounted for in \eq{\PR} due to the neglected uncorrelated component. \par

A pragmatic solution to reduce the risk of excessively amplifying the small spatial scales is to `play it safe' by decreasing \eq{\Prhoo}, even if observational evidence (\textit{e.g.}, from \cite{Desroziers_2005} diagnostics) suggests otherwise, and by inflating the variances by a modest amount. For example, if \eq{\Prhoo = 400}~km instead of 450~km, and \eq{\Psigmao} is inflated by a factor of 2 (smaller than the factor of 4 used with the  diagonal \eq{\PR} in Figure~\ref{fig_convergence_L450}(b)), then the analysis retains a fair share of the improvements seen with an accurate non-diagonal \eq{\PR} while being resilient to perturbations at small spatial scales (Figure~\ref{fig_convergence_unstab}(b)).
A clear drawback with this approach is that it requires empirical tuning and is visibly sub-optimal when comparing Figure~\ref{fig_convergence_unstab}(b) and Figure~\ref{fig_convergence_L450}(c).
A better approach would be to cap the minimum allowed value of the observation-error variances
at the small spatial scales while leaving the large spatial scales unaltered.
At the very small spatial scales, the background and observations are already close to the true state,
so that an inaccurate specification of the observation-error variance at those scales would not have a strong impact on the analysis. For example, if \eq{\sigmao^2\lambok} is capped at \eq{10^{-2}} in Figure~\ref{fig_scale_decomp_L450}(a) then the impact on the analysis error is minor.  \par

Similar ideas have been explored with inter-channel observation-error correlations
\citep{Bormann_2016,Geer_2019,Weston_2014}, by controlling the smallest eigenvalues of an explicit eigen-decomposition of the inter-channel correlation matrix. With other modelling approaches where \eq{\PR} is built explicitly (\textit{e.g.}, from a covariance function) and then inverted directly, another possibility would be to `recondition' \eq{\PR} by adding a scalar multiple of the identity matrix to the matrix modelling the correlated error \mbox{\citep{Tabeart_2020}}. Although this method also alters eigenvalues at large spatial scales, its effect is only significant at small spatial scales when the multiplicative factor is small enough. Producing the same effect with the diffusion model, while retaining a convenient form of \eq{\PR^{-1}}, is more challenging.
 \par

%% file: conclusion.tex
The purpose of this article has been threefold. First, it provides additional insight into the impact of correlated observation error in data assimilation. Second, it addresses the specific problem of accounting for correlated error in along-track (nadir) altimeter observations in a variational data assimilation system (NEMOVAR) for the global ocean. Third, it provides the first practical demonstration of the implicit diffusion approach for modelling observation-error correlations \citep{Guillet_2019} in a realistic data assimilation system.

A key aspect of our study has been to analyse the general data assimilation problem in a simplified context using a Fourier analysis. In doing so, we have highlighted the scale-dependent characteristics of the problem and used them to interpret the results from the simulated data experiments with the ocean data assimilation system. The two key parameters in the implicit diffusion-based correlation model are the length scale \eq{\rho} and the number of diffusion iterations \eq{m}. These parameters control the shape of the correlation function, in particular its spatial range (\eq{\rho}) and the steepness of its spectral slope at small spatial scales (\eq{m}). In designing our experiments, we focussed on studying the sensitivity of the analysis errors to the specification of the basic parameters (\eq{\rhoo}, \eq{\mo}) in the diffusion-based correlation model of the observation-error covariance matrix (\eq{\bR}). The corresponding parameters (\eq{\rhob}, \eq{\mb}) used in the diffusion-based correlation model of the background-error covariance matrix (\eq{\bB}) were fixed to values used in a previous operational version of NEMOVAR \citep{Balmaseda_2013}. Specifically, \eq{\rhob} was spatially varying; \eq{\rhoo} was constant; \eq{\mb} was chosen to model an approximate Gaussian function; and, in most experiments,  \eq{\mo} was chosen to model a SOAR function.

The main findings from our study are as follows. We have shown that assimilating observations with correlated errors, but ignoring them in \eq{\bR} in the data assimilation method, results in overfitting the observations at large spatial scales and underfitting them at small spatial scales. Variance inflation can reduce the overfitting at large spatial scales and is effective when \eq{\rhoo < \rhob}. Nevertheless, as expected, explicitly accounting for the correlations in \eq{\bR} gives best results in terms of reducing the overall analysis error. The improvement is particularly noticeable in zonal velocity, which, via geostrophic balance, benefits from meridional derivative information in the along-track altimeter observations. This information about the small spatial scales is effectively extracted by the observation-error correlation model. When \eq{\rhoo > \rhob}, the derivative information in the observations becomes more important and the impact of the observation-error correlation model is even more impressive at reducing the analysis errors, particularly in zonal velocity.  Variance inflation is shown to be much less effective for the case \eq{\rhoo > \rhob} than the case \eq{\rhoo < \rhob}, as it can only prevent overfitting of observations at large spatial scales by simultaneously hindering the extraction of information at small spatial scales.
 
When \eq{\rhoo < \rhob} (taken here to be roughly 70\% of the minimum value of \eq{\rhob}), the minimization reaches full convergence in half the number of iterations compared to the experiment with a non-inflated diagonal \eq{\bR}. In this case, the non-diagonal \eq{\bR} decreases the condition number. The convergence rate and condition number are  similar to those obtained in the experiment with an optimally inflated diagonal \eq{\bR}.  When \eq{\rhoo > \rhob} (taken here to be roughly twice the \textit{minimum} and slightly larger than the \textit{maximum} value of \eq{\rhob}), the convergence rate with a non-diagonal \eq{\bR} is similar to the one obtained with a non-inflated diagonal \eq{\bR}, but worse than with an optimally inflated \eq{\bR}. Nevertheless, our experiments show that the analysis error is always smaller with a non-diagonal \eq{\bR} than with the optimally inflated diagonal \eq{\bR} at any iteration of the minimisation. However, increasing \eq{\rhoo} even further (to roughly twice the \textit{maximum} value of \eq{\rhob}) significantly increases the condition number and slows down convergence. 

When \eq{\rhoo} is large, the condition number of \eq{\bR} itself becomes large because of small eigenvalues associated with the small spatial scales. In this case, application of \eq{\bR^{-1}} can potentially amplify small-scale features to extreme levels. The risk of this occurring is significantly higher with large values of \eq{\mo} as demonstrated analytically by \mbox{\cite{Goux_2024}} and confirmed by experiments in our study.
For example, when \eq{\mo= 10} (quasi-Gaussian), the condition number of \eq{\bR} becomes huge unless \eq{\rhoo} is chosen to be extremely small.  This extreme sensitivity to large values of \eq{\mo} was the main motivation for the choice of \eq{\mo = 2} (SOAR) in our experiments.\par

In practice, the \eq{\bR} that is used in the data assimilation system will always be an approximation of the true \eq{\bR}. We showed that the accuracy of the analysis is particularly sensitive to errors at small spatial scales that are not fully accounted for in \eq{\bR}.  In this situation,  we found that it was sufficient to reduce \eq{\rhoo} slightly and apply a modest amount of variance inflation to retain much of the benefit to the analysis from using a non-diagonal \eq{\bR} while rendering \eq{\bR} more robust. A better solution would be to modify the diffusion model in a way that would allow more control over the small eigenvalues in \eq{\bR}.

%% file: app_implementation.tex
\begin{figure}
    \centering
    \includegraphics[width=0.5\linewidth]{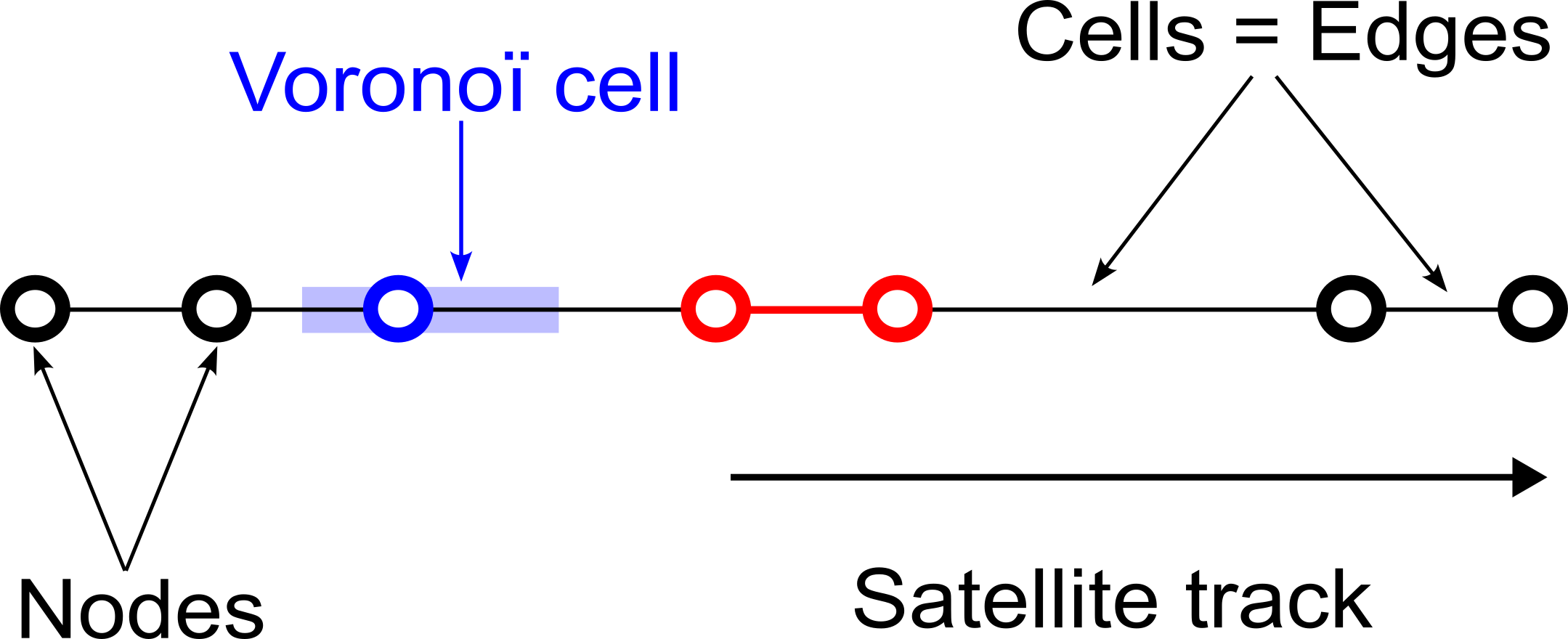}
    \caption{Schematic representation of the mesh used with the 1D diffusion operator for nadir altimeter data. The black circles denote the nodes associated with observations locations. The black lines denote both the edges of the mesh, which connect the observation locations, and the cells of the mesh, which divide the domain into sub-domains. The blue filled-in area is the Voronoï cell associated with the observation location highlighted in blue.}
    \label{fig:mesh_schematic}
\end{figure}

 To discretise the diffusion operator with FEM, we need a mesh that connects the observation locations as illustrated schematically in Figure~\ref{fig:mesh_schematic}. In 1D, the mesh is composed of:
 \begin{itemize}
     \item a set of \eq{N_n} nodes, each corresponding to an observation location;
     \item a set of \eq{N_n-1} edges, each linking two nodes together;
     \item a set of \eq{N_n-1} cells, which partition the domain into sub-domains. 
 \end{itemize}
Edges and cells perfectly overlap in 1D\footnote{n higher dimensions, cells and edges are distinct: in 2D, edges are still segments linking observation locations,  while cells can be triangular surfaces.}. The Voronoï cell, which is not defined explicitly on the mesh, is associated with a node and is delimited by the midpoints of the edges attached to the node (Figure~\ref{fig:mesh_schematic}).\par

Numerically, we can represent the mesh by a connectivity table that indicates which edges are related to which nodes. It can take the form of an array where each row corresponds to an edge, and contains the indices of the two nodes at both its ends. For example, the row associated with the red edge in Figure~\ref{fig:mesh_schematic} would indicate the indices of the nodes highlighted in red. In 1D, building the connectivity table is straightforward. In our case, an edge of the mesh is created between nodes associated with successive observation times. Creating a mesh for 2D observations is more complex, and generally requires a meshing tool such as Atlas\footnote{https://sites.ecmwf.int/docs/atlas/} to create connectivity tables from observation coordinates.

\begin{figure}
    \centering
    \includegraphics[width=0.5\linewidth]{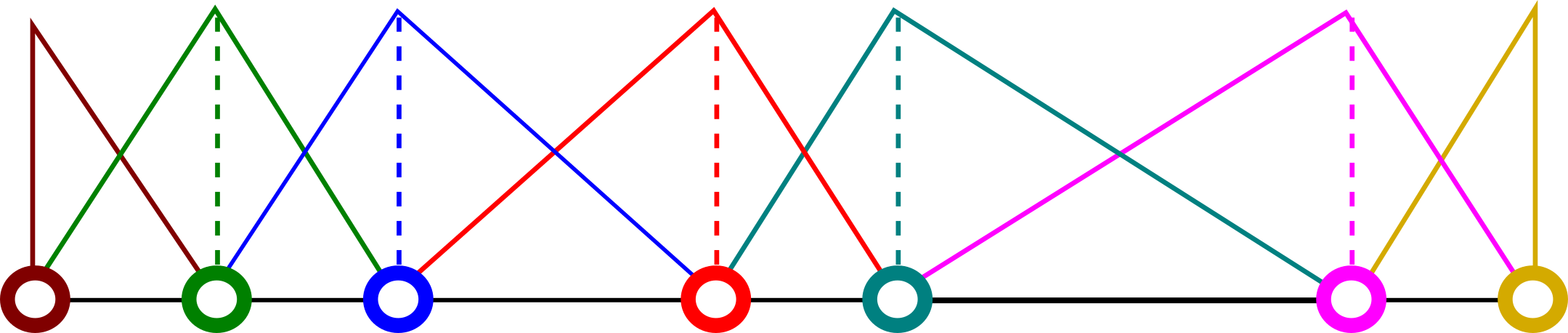}
    \caption{Schematic representation of the shape functions associated in 1D with \eq{\mathcal{P}_1} elements.}
    \label{fig:mesh_P1}
\end{figure}

As defined in Equation~\eqref{eq_D_alt2}, the FEM-based diffusion operator is defined in terms of the mass and stiffness matrices denoted \eq{\bM} and \eq{\bA}, respectively. The elements on row \eq{k} and column \eq{l}  link the \eq{k}-th and \eq{l}-th nodes are defined as \citep[cf. Equations~27 and 28 in][]{Guillet_2019}
\begin{align}
   \bM_{k,l} &= \int_{\Omega} \psi_{k}(z) \psi_{l}(z) dz,\\
   \bA_{k,l} &= \int_{\Omega} \left[\kappa(z) \nabla \psi_{k}(z)\right]^{\rm T}\bm{\nabla}\psi_{l}(\bm{z}) d\bm{z}, \label{eq_MA_integrals}
\end{align}
where \eq{(\psi_{k})_{k\in\llbracket1,N_n\rrbracket}} are the shape functions associated with \eq{\mathcal{P}_1} finite elements, \eq{\kappa(z)} is the diffusion tensor (possibly spatially varying), and \eq{\Omega} is the domain. In the 1D case considered in this study, \eq{\Omega} is the satellite track, \eq{\nabla \psi_{k}} has a single component \eq{\partial \psi_{k} /\partial z}, and \eq{\kappa(z)} is a constant diffusion coefficient.
The function \eq{\psi_k} is piecewise linear on each cell, equal to 1 at the \eq{k}-th node, and 0 at other nodes, as illustrated in Figure~\ref{fig:mesh_P1}. The local compact support of these functions is responsible for the sparse structure of \eq{\bM} and \eq{\bA}: the integrals in Equation~\eqref{eq_MA_integrals} are non-zero only if the nodes \eq{k} and \eq{l} are linked by an edge, or if \eq{k=l}. If the nodes are labelled consistently with their position in the mesh then \eq{\bM} and \eq{\bA} are tri-diagonal matrices. \par

To reduce computational costs, we adopt the common approach of mass lumping \citep[\textit{e.g.}, Section~17.2.4 in][]{Zienkiewicz_2000}, which approximates \eq{\bM} by a diagonal matrix \eq{\bMl}. Each diagonal element of \eq{\bMl} is obtained by summing all the elements in the corresponding row:
\begin{equation}
   \bMl_{k,k} = \sum_{l=1}^{N_n} \bM_{k,l}.
\end{equation}
This assumption is equivalent to taking \eq{\bMl_{k,k}} to be the surface of the Vorono\"{i} cell centred around the \eq{k}-th node (Figure~\ref{fig:mesh_schematic}).\par

The most computationally efficient way to compute the value of the non-zero integrals is to cycle through the cells one by one in order to compute their contribution to all of the integrals. We thus separate the contribution of each cell to the integrals; \textit{i.e.}, we decompose the satellite ground track into the segments  that connect the observation locations:
\begin{align}
   \bM_{n_k,n_l} &= \sum_{\omega\subset\Omega}  \bM_{n_k,n_l}(\omega) \quad \text{with}  \quad\bM_{n_k,n_l}(\omega) =\int_{\omega}\psi_{n_k}(z) \psi_{n_l}(z) dz,
   \label{eq_decomp_M}\\
   \bA_{n_k,n_l} &= \sum_{\omega\subset\Omega}  \bA_{n_k,n_l}(\omega) \quad\; \text{with}  \quad \bA_{n_k,n_l}(\omega) =  \int_{\omega}\left[\kappa(z) \nabla \psi_{n_k}(z)\right]^{\rm T} \nabla\psi_{n_l}(z) dz, \label{eq_decomp_A}
\end{align}
where \eq{\omega} is the sub-domain represented by a cell. As the shape functions are linear within each cell, their gradient is constant on each cell. Using this property, we can simplify the integral in Equation~\eqref{eq_decomp_A} as
\begin{equation}
    \bA_{k,l}(\omega) = \left[ \bm{\nabla} \psi_{k}(\omega) \right]^{\rm T} \left( \int_{\omega}\kappa(\bm{z})  d\bm{z} \right) \bm{\nabla}\psi_{l}(\omega).
    \label{eq_Akl}
\end{equation}
Assuming  \eq{\kappa(z)} is linear on each element (\textit{i.e.}, that it can be written in terms of the basis functions \eq{(\psi_{n_k})_{k\in\llbracket1,N_n\rrbracket}}), the integral in Equation~\eqref{eq_decomp_A} becomes
\begin{equation}
    \int_{\omega}\kappa(\bm{z})  d\bm{z} =  m(\omega)\kappa(\omega),
\end{equation}
where \eq{\kappa(\omega)} is the mean of the values of \eq{\kappa} at the nodes delimiting the cell, and \eq{m(\omega)} is the length of the cell. From Equation~\eqref{eq_Akl}, we have
\begin{equation}
    \bA_{k,l}(\omega) =  m(\omega) \left[ \bm{\nabla}  \psi_{k}(\omega) \right]^{\rm T} \; \kappa(\omega) \;  \bm{\nabla}\psi_{l}(\omega). \label{eq_A_simplex}
\end{equation}

To evaluate the contribution of a cell \eq{\omega} to these integrals, we use the connectivity table to identify the indices \eq{k_1} and \eq{k_2} of the nodes delimiting the cell. Let \eq{z_1} and \eq{z_2} denote the coordinates of these nodes, respectively. For nadir altimeters, observation error is correlated with respect to the observation time along the satellite track. However, we prefer to use the geographical position of the observations to define the coordinate \eq{z} in order to facilitate comparisons with background-error correlations. The observation time is thus converted into a position along the satellite track (with respect to any reference as only relative distances, not absolute positions, are required) by multiplying it by the average ground-scanning velocity of the satellite, provided by the satellite operator.\par

The shape function \eq{\psi_{k_1}} (\eq{\psi_{k_2}}) is equal to 1 at \eq{z_1} (\eq{z_2}) and 0 at \eq{z_2} (\eq{z_1}), and varies linearly between them:
\begin{equation}
    \psi_{k_1}(z) =
\frac{z_{k_2} - z}{m(\omega)} \quad \mbox{and} \quad%
    \psi_{k_2}(z) =  \frac{z - z_{1}}{m(\omega)}, \label{eq_phi_1d}
\end{equation}
where  \eq{m(\omega) = z_{k_2} - z_{k_1}}. The integrals in Equation~\eqref{eq_decomp_M} can be evaluated using Equation~\eqref{eq_phi_1d}. Instead of building a non-diagonal \eq{\bM} and then summing its elements through mass lumping, we build \eq{\bMl} directly:
\begin{equation}
   \left.
 \begin{array}{rcl}
\displaystyle
\bMl_{k_1,k_1}(\omega)
& \displaystyle =  & \displaystyle
 \int_{\omega}\psi_{k_1}(\bm{z}) \psi_{k_1}(\bm{z}) d\bm{z} + \int_{\omega}\psi_{k_1}(\bm{z}) \psi_{k_2}(\bm{z}) d\bm{z} = \frac{m(\omega)}{3} + \frac{m(\omega)}{6} = \displaystyle
\frac{m(\omega)}{2}, \\%
\displaystyle
   \bMl_{k_2,k_2}(\omega)
& \displaystyle =  & \displaystyle
   \int_{\omega}\psi_{k_2}(\bm{z}) \psi_{k_2}(\bm{z}) d\bm{z} + \int_{\omega}\psi_{k_2}(\bm{z}) \psi_{k_1}(\bm{z}) d\bm{z} = \frac{m(\omega)}{3} + \frac{m(\omega)}{6} = \frac{m(\omega)}{2},\\
\displaystyle
   \bMl_{k_1,k_2}(\omega)
& \displaystyle =  & \displaystyle
   0,\\
\displaystyle
   \bMl_{k_2,k_1}(\omega)
& \displaystyle =  & \displaystyle
   0.\label{eq_incr_M}
  \end{array}
   \right\}
\end{equation}

The two components of the gradient of the shape function at nodes \eq{k_1} and \eq{k_2} is
\begin{equation}
    \frac{\partial \psi_{k_{1}}}{\partial z}
    = -\frac{1}{m(\omega)} \quad \mbox{and} \quad
\frac{\partial \psi_{k_2}}{\partial z}
    = \frac{1}{m(\omega)}.
\end{equation}
The integrals in Equation~\eqref{eq_decomp_A} for the stiffness matrix \eq{\bA} then become
\begin{equation}
   \bA_{n_1,n_1}(\omega) = \bA_{n_2,n_2}(\omega) = \frac{\kappa(\omega)}{m(\omega)} \quad \mbox{and} \quad%
   \bA_{n_1,n_2}(\omega) = \bA_{n_2,n_1}(\omega) = -\frac{\kappa(\omega)}{m(\omega)}. \label{eq_incr_A}
\end{equation}

The matrices \eq{\bM} and \eq{\bA} could be built explicitly by initializing their elements to 0, and then, for each cell \eq{\omega}, computing \eq{m(\omega)} and \eq{\kappa(\omega)}, and incrementing the relevant matrix coefficients of \eq{\bMl} and \eq{\bA} based on Equations~\eqref{eq_incr_M} and \eqref{eq_incr_A}. In practice, \eq{\bMl} and \eq{\bA} are not stored as explicit matrices to avoid excessive memory cost. As \eq{\bMl} is diagonal with mass lumping, it can be represented by a single value per node. These values form what is referred to as the \textit{mass field}, which on the \eq{k}-th node  is given by
\begin{equation}
    m(k) = \bMl_{k,k}.
\end{equation}
The similarity with the notation used for the length of each cell (\eq{m(\omega)}) is deliberate, as the mass field can be seen as its counterpart defined on nodes rather than cells.\par

A typical row of a matrix-vector product \eq{\bA\bm{\eta}} is
\begin{equation}
[\bA\bm{\eta}]_{k} = \bA_{k,k}\bm{\eta}_{k} + \bA_{k,k-1}\bm{\eta}_{k-1} + \bA_{k,k+1}\bm{\eta}_{k+1}.
\end{equation}
As the rows of \eq{\bA} sum to 0 by construction, \eq{\bA_{k,k} + \bA_{k,k-1} + \bA_{k,k+1}=0},  we have
\begin{equation}
[\bA\bm{\eta}]_{k} = \bA_{k,k-1}(\bm{\eta}_{k-1} - \bm{\eta_k}) + \bA_{k,k+1}(\bm{\eta}_{k+1} - \bm{\eta_k}), \label{eq_transform_A}
\end{equation}
and there is thus no
need to store the diagonal elements \eq{\bA_{k,k}} of \eq{\bA}. Furthermore, as \eq{\bA} is symmetric, it is sufficient to store only one of its sub-diagonals. Each element on the  sub-diagonal vector can be associated with an edge of the mesh. Instead of building the matrix \eq{\bA}, we build a \textit{stiffness field} defined on the edge \eq{e} associated with the cell \eq{\omega} as
\begin{equation}
    a(e)= - \bA_{k_1,k_2}.
    \label{eq_def_stiffness}
\end{equation}

Knowing the mass field (at nodes) and the stiffness field (on edges) is sufficient to be able to multiply a vector \eq{\bm{\eta}} by \eq{\bI + \bMl^{-1/2}\bA\bMl^{-1/2}} as summarised  in Algorithm~\ref{alg_prod}. This product is the basic operation applied on each (inverse diffusion) iteration in \eq{\bD^{-1}} (Equation~\eqref{eq_D_alt2}). An advantage of this formulation is that the numerical code associated with Algorithm~\ref{alg_prod} for the matrix-vector product is completely generic between 1D and higher dimensions. On the other hand, although the main steps remain the same, the procedure for constructing the mass and stiffness fields depends on the dimension of the domain.

\begin{algorithm}
\caption{Product of a vector \eq{\bm{\eta}} with \eq{\bI + \bMl^{-1/2}\bA\bMl^{-1/2}}}
\label{alg_prod}
\begin{algorithmic}
\State
\State Update the value of \eq{\bm{\eta}} in the halo.
\State
\For {each edge \eq{e_i} with \eq{i \in \llbracket 1,N_n-1 \rrbracket}}
    \State Extract from the connectivity table the indices \eq{k_1}, \eq{k_2} of the nodes at both ends of the edge \eq{e_i}
    \State  $f(e_i) \gets a(e_i) \left ( \frac{\eta_{k_1}}{\sqrt{m_{k_1}}} - \frac{\eta_{k_2}}{\sqrt{m_{k_2}}}\right);$\Comment{Compute the flux through the edge}
\EndFor
\State
\For {each edge \eq{e_i} with \eq{i \in \llbracket 1,N_e \rrbracket}}
    \State $\eta_{k_1} \gets \eta_{k_1} - \frac{f(e_i)}{\sqrt{m_{k_1}}};$ \Comment{Apply the flux to the first node}
    \State $\eta_{k_2} \gets \eta_{k_2} + \frac{f(e_i)}{\sqrt{m_{k_2}}};$ \Comment{Apply the flux to the second node}
\EndFor
\end{algorithmic}
\end{algorithm}